\documentclass[letterpaper]{article}
\usepackage{amssymb, amsmath,tikz-cd,mathrsfs,amsthm}
\usepackage{mathtools}
\usepackage{comment}
\usepackage{cite}
\usepackage{hyperref}

\hypersetup{
	unicode=false,          
	pdftoolbar=true,        
	pdfmenubar=true,        
	pdffitwindow=false,     
	pdfstartview={FitH},    
	pdftitle={Hodge Numbers of O'Grady 6 Via Ng\^o Strings},    
	pdfauthor={Ben Wu},     
	pdfnewwindow=true,      
	colorlinks=false,       
	linkcolor=red,          
	citecolor=green,        
	filecolor=cyan,         
	urlcolor=magenta        
}

\newcommand{\bbC}{\mathbb{C}}
\newcommand{\bbP}{\mathbb{P}}
\newcommand{\bbQ}{\mathbb{Q}}

\newcommand{\bbZ}{\mathbb{Z}}

\newcommand{\bfa}{\mathbf{a}}

\newcommand{\bfG}{\mathbf{G}}
\newcommand{\bfM}{\mathbf{M}}

\newcommand{\scrA}{\mathscr{A}}
\newcommand{\scrI}{\mathscr{I}}
\newcommand{\scrL}{\mathscr{L}}

\newcommand{\scrF}{\mathscr{F}}
\newcommand{\scrG}{\mathscr{G}}
\newcommand{\scrIC}{\mathscr{IC}}

\newcommand{\tildea}{\widetilde{a}}
\newcommand{\tildei}{\widetilde{i}}
\newcommand{\tildem}{\widetilde{m}}
\newcommand{\tildeC}{\widetilde{C}}
\newcommand{\tildecalC}{\widetilde{\mathcal{C}}}
\newcommand{\tildeD}{\widetilde{D}}
\newcommand{\tildecalD}{\widetilde{\mathcal{D}}}

\newcommand{\tildeM}{\widetilde{M}}

\newcommand{\tildephi}{\widetilde{\phi}}
\newcommand{\tildesigma}{\widetilde{\sigma}}
\newcommand{\Jhat}{\widehat{J}}

\newcommand{\tildeJ}{\widetilde{J}}
\newcommand{\tildeK}{\widetilde{K}}

\newcommand{\tildebfa}{\widetilde{\bf{a}}}

\newcommand{\tildeiota}{\widetilde{\iota}}
\newcommand{\tildealb}{\widetilde{\textup{alb}}}
\newcommand{\rat}{\textup{rat}}
\DeclareMathOperator{\codim}{codim}
\DeclareMathOperator{\Sym}{Sym}

\newcommand{\et}{\textup{et}}
\newcommand{\red}{\textup{red}}

\newcommand{\calA}{\mathcal{{A}}}
\newcommand{\calC}{\mathcal{C}}

\newcommand{\calH}{\mathcal{H}}
\newcommand{\calI}{\mathcal{I}}

\newcommand{\calO}{\mathcal{O}}

\newcommand{\fraks}{\mathfrak{s}}

\newcommand{\frakS}{\mathfrak{S}}
\newcommand{\frakM}{\mathfrak{M}}

\newcommand{\Npq}{N_{pq}}

\newcommand{\Upq}{U_{pq}}
\newcommand{\Vpq}{V_{pq}}

\newcommand{\Aut}{\textup{Aut}}
\newcommand{\Tr}{\textup{Tr}}
\newcommand{\Irr}{\textup{Irr}}
\newcommand{\Obj}{\textup{Obj}}

\newcommand{\inv}{^{-1}}

\newcommand{\Ext}{\textup{Ext}}
\newcommand{\Fix}{\textup{Fix}}
\newcommand{\Pic}{\textup{Pic}}

\newcommand{\Cred}{C_{\red}}

\newcommand{\sheafEnd}{\mathcal{E}nd}

\newcommand{\leqp}{%
	\mathrel{\raisebox{-0.5ex}{$\scriptscriptstyle($}}%
	\leq
	\mathrel{\raisebox{-0.5ex}{$\scriptscriptstyle)$}}%
}

\newtheorem{thm}{Theorem}
\newtheorem{defn}{Definition}[subsection]
\newtheorem{lem}[defn]{Lemma}
\newtheorem{prop}[defn]{Proposition}
\newtheorem{cor}[defn]{Corollary}
\newtheorem{rmk}[defn]{Remark}

\newtheorem*{fact}{Fact}
\newtheorem*{NST}{Ng\^o Support Theorem}

\title{Hodge Numbers of O'Grady 6 Via Ng\^o Strings}
\author{Ben Wu}
\date{}

\begin{document}
	
	\maketitle
	
	\begin{abstract}
		We give an alternative computation of the Betti and Hodge numbers for manifolds of $OG6$ type using the method of Ng\^o Strings introduced by de Cataldo, Rapagnetta, and Sacc\`a. 
	\end{abstract}
	
	\tableofcontents

\section{Introduction}
The Betti and Hodge numbers of manifolds of $OG6$ type were first computed by Mongardi, Rapagnetta, and Sacc\`a in \cite{MRS1}. The purpose of this paper is to provide a new computation of the Betti and Hodge numbers of these manifolds using the method of Ng\^o strings introduced by de Cataldo, Rapagnetta, and Sacc\`a in \cite{dCRS}. We then prove Theorem \ref{main thm hodge structure}, which gives a description of the Hodge structure of some special manifolds of $OG6$ type in terms of the Hodge structure of an Abelian surface. Although we do not describe the Hodge structure of all manifolds of $OG6$ type, as is done in \cite{Green-Kim-Laza-Robles}, we note that our computation of the Hodge structure using Ng\^o strings does not require prior knowledge of the Hodge numbers as an input.

The computation of the Betti and Hodge numbers follows the strategy in \cite{dCRS} for the $OG10$ case and is summarized below. We consider two special manifolds, both admitting Lagrangian fibrations onto $\bbP^3$. The first, denoted by $\tildeM$, is of $OG6$ type and the second, denoted by $N$, is deformation equivalent to a generalized Kummer variety of dimension six. We note that the Lagrangian fibration $\tildeM \to \bbP^3$ is the same as the ones considered in \cite{R1, MRS1}. The two fibrations $\tildeM \to \bbP^3$ and $N \to \bbP^3$ are linked via a group scheme $G \to \bbP^3$ which acts fiberwise on both $\tildeM$ and $N$. In particular, $\tildeM$ contains an open dense subvariety which is a $G$-torsor over the entire base and $N$ contains an open dense subvariety which is a $G$-torsor over an open dense subvariety strictly contained in the base.  Given this relationship between $\tildeM$ and $N$, we are able to use the known cohomology of generalized Kummer varieties, as well as some smaller dimensional fibrations, to deduce the Betti and Hodge numbers of $\tildeM$ by working in an appropriate Grothendieck group. Noting that the pure Hodge structure of $N$ and these smaller dimensional fibrations can all be described in terms of the Hodge structure of an Abelian surface $J$, we then describe the Hodge structure of $\tildeM$ in terms of the Hodge structure of $J$ (see Theorem \ref{main thm hodge structure}). As in the $OG10$ case, the main technical tool that allows us to compare the cohomology of $\tildeM$ and $N$ is the Ng\^o Support Theorem (see Section \ref{Ngo support thm section}). A key difference with the $OG10$ case is that the fibers of the group scheme $G$ become disconnected over various loci in the base.

\section*{Acknowledgments}
I would first like to greatly thank my advisor Mark de Cataldo for the many enlightening discussions and his continued support throughout the writing of this paper. I would also like to thank Mads Bach Villadsen for the many useful discussions and for reading several preliminary drafts, Lisa Marquand for the many useful discussions and her support, and Yoonjoo Kim for a useful discussion on the $LLV$ decomposition. I would finally like to give a special thanks to Antonio Rapagnetta for an extremely useful discussion on monodromy (in particular, see Lemma \ref{O(p1-p2) not in identity component cor} and Section \ref{monodromy of irred comp section}).

\section{Preliminaries}

\subsection{Notation}
We work over the field of complex numbers. A variety is a separated scheme of finite type over $\bbC$. By point, we always mean closed point. Given a morphism of varieties $f \colon X \to Y$ and a subvariety $Z \subseteq Y$, we set $X_Z \coloneqq f\inv(Z)$. In particular, if $Z=\{y\}$ is a point, then $X_y = f\inv(y)$ denotes the fiber. 

Given a variety $X$, we work with the bounded constructible derived category, $D^b(X,\bbQ)$ (see \cite{dCM1} for the basics and references) and Saito's bounded derived category of algebraic mixed Hodge modules which rational structure, $D^bMHM_{alg}(X)$ (see \cite{Sai-1990} for the foundations and \cite{Sch-2014} for basics and references). There is an exact functor $\rat \colon D^bMHM_{alg}(X) \to D^b(X, \bbQ)$ which takes the standard $t$-structure on $D^bMHM_{alg}(X)$ to the middle perversity $t$-structure on $D^b(X, \bbQ)$. An important fact is that a splitting of an object $K \in D^bMHM_{alg}(X)$ induces a splitting of the object $\rat(K) \in D^b(X,\bbQ)$.

If $Z \subseteq X$ is a closed irreducible subvariety and $\scrL$ is a polarizable variation of rational pure Hodge structures of weight $w(\scrL)$ on some Zariski dense open subset $V \subseteq Z^{reg}$, then the intersection cohomology object $IC_{Z}(\scrL) \in MHM_{alg}(X)$ yields, via $\rat$, the usual perverse intersection cohomology complex of $Z$ with coefficients in the underlying local system $\scrL$. We denote by $\scrIC_Z(\scrL) \coloneqq IC_Z(\scrL)[-\dim(Z)]$ the non-perverse intersection complex; this is useful when making geometric statements.

\subsection{The Ng\^o Support Theorem}\label{Ngo support thm section}
We begin by defining the notion of a $\delta$-regular weak Abelian fibration, first introduced by B.C. Ng\^o in \cite{Ngo-2010}. Let $g \colon G \to B$ be a smooth commutative group scheme over $B$ and $g^0 \colon G^0 \to B$ be the identity component. Given a point $b \in B$, there is the canonical Chevalley devissage short exact sequence
\begin{equation*}\label{chevalley devissage of fibers ses}
0 \to R_b \to G^0_b \to A_b \to 0
\end{equation*}
of the fiber $G^0_b$ of $G^0$ over the perfect field $k(b)$, where $A_b$ is an Abelian variety and $R_b$ is affine, connected, and maximal with respect to these properties. 

The function $\delta \colon B \to \bbZ$ sending $b \in B$ to $\dim_{k(b)} R_b$
is upper-semicontinuous and we define the $\delta$-loci $B_i$ by $B_i \coloneqq \{b \in B \mid \delta(b) = i\}$. Now suppose that $d = \dim_{k(b)} G_b$ is constant on $B$ (it is so on the connected components of $B$ by smoothness). For any fixed prime $\ell$, the Tate module of $G/B$, denoted by $T_{\et, \overline{\bbQ}_\ell}(G/B)$, is defined to be the $\overline{\bbQ}_\ell$-adic counterpart of the object $R^{2d-1}g_!^0 \bbQ_{G^0}(d)$, i.e. it is defined by the same formula using the \'etale topology/cohomology formalism. The Tate module $T_{\et, \overline{\bbQ}_\ell}(G/B)$ is defined to be polarizable if \'etale locally, there is a pairing
\begin{equation*}
T_{\et, \overline{\bbQ}_\ell}(G/B) \otimes T_{\et, \overline{\bbQ}_\ell}(G/B) \to \overline{\bbQ}_\ell(1)
\end{equation*}
such that for every $b \in B$, the kernel of the pairing at $b$ is $T_{\et, \overline{\bbQ}_\ell}(R_b)$. 

\begin{defn}
	A weak Abelian fibration, denoted as a triple $(M,B,G)$ is a pair of of morphisms $M \xrightarrow{f} B \xleftarrow{g}G$ such that $(a)$ $f$ is proper; $(b)$ $G$ is a smooth commutative group scheme over $B$; $(c)$ $f$ and $g$ have the same pure relative dimension $d$; $(d)$ there is an action $a \colon G \times_B M \to M$ of $G$ on $M$ over $B$; $(e)$ the action has affine stabilizers at every point $m \in M$; $(f)$ the Tate module $T_{\et, \overline{\bbQ}_\ell}(G/B)$ is polarizable. If in addition the $\delta$-loci, $B_i$, satisfy the inequality $\codim B_i \geq i$ for every non-negative integer $i$, then $(M,B,G)$ is said to be a $\delta$-regular weak Abelian fibration. In particular, this inequality implies that the general fiber $G^0_b$ is an Abelian variety.
\end{defn}

We now state the refinement of the Ng\^o Support Theorem given by de Cataldo, Rapagnetta, and Sacc\`a in \cite{dCRS}. We work with the following set up. Let $(M,B,G)$ be a $\delta$-regular weak Abelian fibration. Suppose that $d = \dim(M/B) = \dim(G/B)$, $M/B$ is projective, $B$ is irreducible, and $M$ is a rational homology manifold, i.e. $\scrIC_M \simeq \bbQ_M$.

Let $\scrA$ be a finite set enumerating the supports of $Rf_*\bbQ_M$, so that the supports are denoted $Z_{\alpha}$ for $\alpha \in \scrA$. There is an open dense subvariety $V_{\alpha} \subset Z_{\alpha}$ and a short exact sequence of smooth commutative algebraic group schemes over $V_{\alpha}$ with connected fibers
\begin{equation*}
0 \to R_{V_\alpha} \to G^0|_{V_{\alpha}} \to A_{V_{\alpha}} \to 0
\end{equation*}
realizing the Chevalley devissage point-by-point on $V_{\alpha}$ (cf. \cite[\S 7.4.8]{Ngo-2010}). Denote the natural map $A_{V_\alpha} \to V_\alpha$ by $g_\alpha$ and let $\delta^{ab}_\alpha$ denote the relative dimension of $A_{V_\alpha}/V_\alpha$. Let $\Lambda_\alpha^\bullet \coloneqq R^\bullet g_{\alpha*} \bbQ_{A_{V_\alpha}}$ where $0 \leq \bullet \leq \delta^{ab}_{\alpha}.$

\begin{NST}
	Let $\langle \bullet \rangle \coloneqq [-2\bullet](-\bullet)$. With the set-up above, there is an isomorphism 
	\begin{equation*}
	Rf_*\bbQ_M \simeq \bigoplus_{\alpha \in \scrA} \scrI_\alpha \langle d-\delta^{ab}_\alpha \rangle
	\end{equation*}
	in $D^bMHM_{alg}(B)$ of pure objects of weight $0$, where 
	\begin{equation*}
	\scrI_\alpha \coloneqq \bigoplus_{\bullet = 0}^{2\delta^{ab}_\alpha} \scrIC_{Z_\alpha}(\Lambda_\alpha^\bullet \otimes \scrL_\alpha)[-\bullet]
	\end{equation*}
	and $\scrL_\alpha$ is a polarizable variation of pure Hodge structures of weight $0$. The summands $\scrI_\alpha$ are called Ng\^o strings. 
\end{NST}

\subsection{The Varieties \texorpdfstring{$\tildeM$}{tildeM}, \texorpdfstring{$M$}{M}, \texorpdfstring{$N$}{N}, \texorpdfstring{$B$}{B} and The Group Scheme \texorpdfstring{$G$}{G}}\label{Kummer and OG type moduli spaces section}
In this section, we discuss the construction of the varieties $\tildeM, M, N, B$ and the group scheme $G$. In the remainder of this paper, $J \coloneqq J(C_0)$ will always be the Jacobian of a fixed genus 2 curve $C_0$ with $NS(J) = \bbZ\theta$, where $\theta$ is a symmetric theta divisor, i.e. $\theta$ is invariant under pullback by the natural $(-1)$-involution on $J$. Using the theta divisor, we can identify $J$ with the dual Abelian variety, $J^\vee$ by sending a point $x \in J$ to $\theta_x - \theta$ where $\theta_x \coloneqq \theta+x$.  

We now recall some basic definitions about pure dimension one sheaves on the principally polarized Abelian surface $(J, \theta)$. A coherent sheaf $\scrF$ on $J$ is pure of dimension one if for every non-trivial subsheaf $\scrG \subset \scrF$, the dimension of the schematic support of $\scrG$ is one. A pure dimension one sheaf $\scrF$ is Gieseker (semi-)stable with respect to $\theta$ if for all proper pure dimension one quotients $\scrF \twoheadrightarrow \scrG$, the following inequality holds,
\begin{equation*}
\frac{\chi(\scrF)}{c_1(\scrF) \cdot \theta} \leqp \frac{\chi(\scrG)}{c_1(\scrG) \cdot \theta}.
\end{equation*}

If $\scrF$ is a pure dimension one sheaf, a key property is that there exists a length one resolution by locally free sheaves of the same rank (cf. \cite[\S 1.1]{HL}). Using this resolution, Le Potier defines in \cite[\S 2.2]{LP} the Fitting support of $\scrF$ as the vanishing subscheme of the induced morphism between determinant bundles. We note that this definition is independent of the resolution. 

\begin{rmk}
	The Fitting support of a pure dimension one sheaf $\scrF$ on $J$ contains the schematic support of $\scrF$ and represents $c_1(\scrF)$ in the Chow ring of $J$. In particular, two pure dimension one sheaves $\scrF$ and $\scrG$ have isomorphic determinant bundles if and only if their Fitting supports are linearly equivalent.
\end{rmk}

Consider a Mukai vector of the form $v = (0, 2\theta, \chi) \in H^{ev}_{alg}(J, \bbZ)$ with $\chi \neq 0$. Sheaves on $J$ with such a Mukai vector are pure of dimension one. Let $\bfM_v$ be the moduli space of semi-stable sheaves on $J$ with Mukai vector $v$. By \cite{Mukai}, $\bfM_v$ is a normal irreducible projective variety of dimension $10$ and the smooth locus, which is precisely the locus parameterizing strictly stable sheaves, admits a symplectic form. There is a natural morphism
\begin{equation}\label{M_v to J^vee times J morphism}
	\bfa_v \colon \bfM_v \to J^\vee \times J; \quad [\scrF] \mapsto \big(\det(\scrF) \otimes \calO_J(2\theta)^\vee, \textstyle \sum c_2(\scrF)\big),
\end{equation}
where $[\scrF]$ is the $S$-equivalence class of $\scrF$ (see \cite[Definition 1.5.3]{HL}), $\det(\scrF)$ is the determinant bundle of $\scrF$, and $\sum c_2(\scrF)$ is the sum in $J$ of any representative of $c_2(\scrF)$ in $CH_0(J)$, the Chow group of zero cycles on $J$ (cf. \cite[\S 4.5 and Proposition 10.3.6 ]{HL}). 

The fiber of this morphism	$M_v \coloneqq \bfa_v\inv(\calO_J, 0)$ will be our main object of interest in this paper. The points of $M_v$ parameterize $S$-equivalence classes of semi-stable sheaves, or equivalently isomorphism classes of polystable sheaves (see \cite[Definitions 1.5.4]{HL}), on $J$ having determinant bundle isomorphic to $\calO_J(2\theta)$ and second Chern class summing up to 0.

The moduli spaces $\bfM_v$ and the fibers $M_v$ depend on the parity of the integer $\chi$ in the Mukai vector $v = (0,2\theta,\chi)$. If $\chi$ is odd, then the Mukai vector $v$ is primitive and  $\bfM_v$ is smooth of dimension 10 and the fiber $M_v$ is deformation equivalent to a generalized Kummer variety of dimension 6 (see \cite[Theorem 0.1 and Theorem 0.2]{Y1}). If $\chi$ is even, then the Mukai vector $v = 2v'$ is twice a primitive Mukai vector and Rapagnetta shows in \cite{R1} that the moduli space $\bfM_v$ and the fiber $M_v$ are both reduced, but singular. The singular locus $\Sigma_v \subset M_v$ parameterizes polystable sheaves of the form $\scrF_1 \oplus \scrF_2$ with $v(\scrF_i) = v'$, determinant bundle isomorphic to $\calO_J(2\theta)$ and second Chern class summing up to 0. The locus $\Sigma_v$ is itself singular and its singular locus, $\Omega_v \subset \Sigma_v$, parameterizes polystable sheaves of the form $\scrF \oplus \scrF$. The fiber $M_v$ admits a symplectic resolution $\pi_v \colon \tildeM_v \to M_v$, where $\tildeM_v$ is an irreducible holomorphic symplectic manifold which is birational, and hence deformation equivalent, to O'Grady's six dimensional exceptional example (see \cite[Proposition 2.2.1]{R1}). 

\begin{rmk}\label{fibers of symplectic resolution}(cf. \cite[Remark 1.1.5]{R1})
	Let $\pi \colon \tildeM_v \to M_v$ be the symplectic resolution. If $[\scrF] \in \Sigma \setminus \Omega \subset M$, then $\pi\inv([\scrF]) \simeq \bbP^1$. If $[\scrF] \in \Omega \subset M$, then $\pi\inv([\scrF])$ is a smooth 3-dimensional quadric. 
\end{rmk}

\begin{prop}\label{DT for symplectic resolution}
	There is a canonical isomorphism in $D^bMHM_{alg}(M)$ and in $D^b_c(M,\bbQ)$ (turn off the Tate twists)
	\begin{equation*}
		R\pi_*\bbQ_{\tildeM} \simeq \scrIC_{M} \oplus \bbQ_{\Sigma}[-2](-1) \oplus \bbQ_\Omega[-6](-3).
	\end{equation*}
	
	\begin{proof}
		By \cite{Kaledin}, the symplectic resolution $\pi$ is semismall. Rapagnetta, in the proof of \cite[Theorem 2.1.7]{R1}, shows that the singular locus $\Sigma \subset M$ is isomorphic to $(J^\vee \times J)/ \pm 1$. In particular, $\Sigma$ has finite quotient singularities, which implies that the intersection complex $\scrIC_{\Sigma}$ is the constant sheaf. The result then follows from the Decomposition Theorem for semismall morphisms (see Theorem 4.2.7 in the survey \cite{dCM1}). 
	\end{proof}
\end{prop}

There is a Le Potier support morphism from the moduli space $\bfM_v$ to the Hilbert scheme parameterizing closed subschemes of $J$ with cohomology class $2\theta$ (cf. \cite[\S 2.2]{LP}). Since all sheaves parameterized by $M_v$ have isomorphic determinant bundles, the image of the support morphism restricted to $M_v \subset \bfM_v$ is the linear system $|2\theta|$. In the remainder of the paper, we fix some Mukai vectors and will work with the following varieties. Set 
\begin{equation}\label{M,N def}
M \coloneqq M_{(0,2\theta, -2)}; \quad N \coloneqq M_{(0,2\theta,-3)}
\end{equation}
and let $\pi \colon \tildeM \to M$ be the symplectic resolution. Set $B \coloneqq |2\theta| \simeq \bbP^3$, let $m \colon M \to B$ and $n \colon N \to B$ denote the respective support morphisms, and let $\tildem \colon \tildeM \xrightarrow{\pi} M \xrightarrow{m} B$ denote the composition.

\begin{rmk} \label{flatness of the Lagrangian fibrations}
	The morphisms $\tildem$ and $n$ are Lagrangian fibrations by \cite[Theorem 1]{Mat2} and are equidimensional by \cite[Theorem 1]{Mat1} and hence flat by \cite[Theorem 23.1]{Matsumura}. Moreover, 
	\begin{enumerate}
		\item since $M$ has a symplectic resolution, it has canonical singularities and is thus Cohen-Macaulay (see \cite[Corollary 5.24 and Lemma 5.12]{KM});
		\item since the fibers of $m$ are dominated by the corresponding fibers of $\tildem$, they all have the same dimension 3.
	\end{enumerate}
	Since the base $B$ is smooth, it again follows again from by \cite[Theorem 23.1]{Matsumura} that the morphism $m \colon M \to B$ is flat.
\end{rmk}

We now discuss the group scheme $G$. Recall from Equation \ref{M_v to J^vee times J morphism} that there is a map $\bfa_{(0,2\theta,-4)} \colon \bfM_{(0,2\theta,-4)} \to J^\vee \times J$ from the moduli space of sheaves with Mukai vector $(0,2\theta,-4)$ to the Abelian fourfold $J^\vee \times J$. Let $M_{(0,2\theta,-4)} \coloneqq \bfa_{(0,2\theta,-4)} \inv(\calO_J,0)$ be the corresponding fiber over $(\calO_J, 0) \in J^\vee \times J$ and let $\bfM_{(0,2\theta,-4)}^B \coloneqq \bfa_{(0,2\theta,-4)}\inv(\{\calO_J\} \times J)$ be the pre-image of the slice $\{\calO_J\} \times J \subset J^\vee \times J$. Note that since the sheaves parameterized by $\bfM_{(0,2\theta,-4)}^B$ have determinant bundles isomorphic to $\calO_J(2\theta)$, the Fitting support of these sheaves all lie in the linear system $B=|2\theta|$. Let $\bfG \subset \bfM_{(0,2\theta,-4)}$ be the locus parameterizing sheaves which are the pushforwards of line bundles from their schematic supports, $\bfG^B \coloneqq  \bfG \cap \bfM_{(0,2\theta,-4)}^B$, 
\begin{equation}\label{group scheme def}
	G \coloneqq \bfG \cap M_{(0,2\theta,-4)},
\end{equation}
and let $g \colon G \to B$ be the restriction of the support morphism to $G$. 

Following the proof of Lemma 3.3.1 in \cite{dCRS}, one can show that $\bfG$ is a nonempty Zariski open subset of $\bfM_{(0,2\theta,-4)}$. Following the proof of Corollary 3.3.3 in \cite{dCRS}, one can show the open subset $\bfG^B \subset \bfM_{(0,2\theta,-4)}^B$ can be identified with the relative degree-$0$ Picard scheme $\Pic^0_{\calC/B}$ of the family $\calC/B$ of curves in the linear system $B$. Using this identification, the morphism $\bfM_{(0,2\theta,-4)} \xrightarrow{\bfa_{(0,2\theta,-4)}} J^\vee \times J \xrightarrow{pr_2} J$ restricted to $\bfG^B$ induces a morphism $\bfa \colon \Pic^0_{\calC/B} \to J \times B$ which fiberwise is the map
\begin{equation}\label{a equation}
	a \colon \Pic^0(C_b) \to J; \quad L \mapsto \sum c_2(i_{b*}L).
\end{equation}
where $C_b$ is the curve in $J$ corresponding to $b \in B$ and $i_b \colon C_b \to J$ is the inclusion. A computation using Riemann-Roch without denominators for regular embeddings (cf. \cite[Theorem 15.3 and Example 15.3.6]{F1}) shows that the map in Equation \ref{a equation} is a group homomorphism. An immediate corollary is the following.

\begin{prop}
	The nonempty open subset $G \subset M_{(0,2\theta,-4)}$ can be identified with the kernel of the map $\bfa \colon \Pic^0_{\calC/B} \to J \times B$ and thus inherits the structure of a $B$-group scheme.
\end{prop} 

\begin{prop}\label{(tildeM,B,G), (N,B,G) are drwaf}
	The triples $(\tildeM,B,G)$ and $(N,B,G)$ are $\delta$-regular weak Abelian fibrations.
	\begin{proof}
		We give the proof for the triple $(\tildeM,B,G)$ noting that the proof for the triple $(N,B,G)$ is analogous. Let $\bfM \coloneqq \bfM_{(0,2\theta,-2)}$. By Lemma 3.5.1 in \cite{dCRS}, tensoring a sheaf by a line bundle with the same Fitting support induces an action $\bfa_{\bfM} \colon \bfG \times_{B} \bfM \to \bf M$,	which restricts to an action $a_{M} \colon G \times_{B} M \to M$. The action $a_M$ lifts to a unique action $\tildea_M \colon G \times_{B} \tildeM \to \tildeM$ since the resolution $\tildeM \to M$ is the blowing up of a locus which is invariant under the action of the group scheme. The proof of Lemma 3.5.4 in \cite{dCRS} can be adapted to our case to show that the action of $G$ on $M, \tildeM,$ and $N$ have affine stabilizers. $\delta$-regularity follows from the description of the group scheme given in the upcoming Propositions \ref{irred components over integral curves prop}, \ref{irred comp of tildeM, N over R(1) u R(2)}, and \ref{irred comp of tildeM, N over NR}.
		 
		Polarizability of the Tate module for $G$ can be proven by adapting de Cataldo's proof of polarizability for the relative Prym variety in \cite[Theorem 4.7.2]{dC1}. The key observation is the following. For any $b \in B$, let $C_b = \sum m_{b_j} C_{b_j}^{\red} \subset J$ be the corresponding curve in $J$ and let $\tildeC_{b_j} \to C_{b_j}^{\red}$ be the normalization of the irreducible components of $C_b$. Let $\tildealb_{b_j} \colon \Pic^0(\tildeC_{b_j}) \to J$ be the morphism on Albanese varieties induced by the map $\tildei_{b_j} \colon \tildeC_{b_j} \to J$ and let $\tildea_{b} \coloneqq \sum m_{b_j} \tildealb_{b_j}$. One then checks that the morphisms $\tildea_b$ and $\tildei^*_b$ play the same role as the morphisms $N^{ab}_{p_a}$ and $\widetilde{p}_a^*$ in Lemma 4.7.1 of \cite{dC1}. One then follows the proof of Theorem 4.7.2 in \cite{dC1} to conclude. 
	\end{proof}
\end{prop}

\section{Irreducible Components}\label{irred components section}
In this section, we first describe Rapagnetta's stratification of the linear system $B = |2\theta|$ and then study the fibers of $\tildem \colon \tildeM \to B$, $m \colon M \to B$ and $n \colon N \to B$ over each stratum. The main results of this section are summarized in Proposition \ref{tildeM,M,N irred components prop}.

\begin{prop}\cite[Proposition 2.1.3]{R1}\label{stratification of |2theta|}
	Rapagnetta's stratification of the linear system $|2\theta|$ according to analytic type of singularity is the following:
	\begin{itemize}
		\item Stratum S: the open dense locus parameterizing smooth curves of genus $5$.
		
		\item Stratum N(1): the locus parameterizing irreducible nodal curves singular in a unique 2-torsion point. $\overline{N(1)} = \cup_{p \in J[2]} N_p$ where $N_p \simeq \bbP^2$ parameterizes curves singular at $p$.
		
		\item Stratum N(2): the locus parameterizing irreducible nodal curves singular in exactly two distinct 2-torsion points. $\overline{N(2)} = \cup_{p,q \in J[2], p\neq q} \Npq$ where $\Npq \coloneqq N_p \cap N_q \simeq \bbP^1$. Note that there are $\binom{16}{2} = 120$ lines $\Npq$.
		
		\item Stratum N(3): the locus parameterizing irreducible nodal curves singular in exactly three distinct 2-torsion points. We have $|N(3)| = 240$.
		
		\item Stratum R(1): the locus parameterizing reducible curves of the form $\theta_x +\theta_{-x}$ with two distinct singular points. 
		
		\item Stratum R(2): the locus parameterizing reducible curves of the form $\theta_x + \theta_{-x}$ with a unique singular point belonging to $J[2]$. 
		
		\item Stratum NR: the locus parameterizing non-reduced curves of the form $C = 2\Cred$ where $\Cred \simeq \theta_p$ is the translate of the theta divisor by a 2-torsion point $p \in J[2]$. We have $|NR| = 16$. 
	\end{itemize}
	In particular, the locus parameterizing singular curves in $B$ consists of $17$ irreducible divisors, namely the divisor $R = \overline{R(1)}$ parameterizing reducible curves and the divisors $N_p$ for $p \in J[2]$ parameterizing curves which are singular at $p$. 
\end{prop}

\begin{rmk}\label{poset structure of |2theta|}
	The following poset structure of Rapagnetta's stratification is implicit in the description given by Rapagnetta in \cite[Proposition 2.1.3]{R1}: 
	\begin{equation*}
		\begin{tikzcd}
			^0NR \arrow[r, hook] \arrow[dr, hook] & ^1\overline{R(2)} \arrow[r, hook] & ^2R \arrow[r, hook] &^3B\\
			^0N(3) \arrow[r, hook] & ^1\overline{N(2)} \arrow[r, hook] & ^2\overline{N(1)} \arrow[ur, hook]
		\end{tikzcd}
	\end{equation*}
	where the superscripts on the left denote the respective dimensions of the strata.
\end{rmk}

\begin{prop}\label{tildeM,M,N irred components prop}
	The number of irreducible components of the fibers of $\tildem \colon \tildeM \to B$, $m \colon M \to B$ and $n \colon N \to B$  over each stratum are summarized by the following table. The superscript on the left of each stratum denotes the respective dimension and the entries in the table are the number of irreducible components in the respective fibers.
	\begin{center}
		\begin{tabular}{ |c|c|c|c|c|c|c|c| }
			\hline
			& &  &  &  &  &  & \\[-1em]
			Strata & $^3S$ & $^2N(1)$ & $^1N(2)$ & $^0N(3$) & $^2R(1)$ & $^1R(2)$ & $^0NR$\\ \hline
			& &  &  &  &  &  & \\[-1em]
			$\tildeM_C$ & 1 & 1 & 2 & 4 & 2 & 2 & 34\\ \hline
			& &  &  &  &  &  & \\[-1em]
			$M_C$ & 1 & 1 & 2 & 4 & 1 & 1 & 17\\ \hline
			& &  &  &  &  &  & \\[-1em]
			$N_C$ & 1 & 1 & 2 & 4 & 2 & 2 & 2\\ 
			\hline 
		\end{tabular}
	\end{center}
\end{prop}

\subsection{Irreducible Components Over Integral Curves}\label{irred components over integral curves section}
Fix an integral curve $C \in B$ with $k$ nodes. Let $i \colon C \to J$ be the inclusion map, $\nu : \tildeC \to C$ be the normalization of $C$, and $a \colon \Pic^0(C) \to J$ be the group homomorphism sending $L \in \Pic^0(C)$ to $\sum c_2(i_*L)$ described in Equation \ref{a equation}. Since $C$ is integral, all sheaves parameterized by the fiber $M_C$ are stable. In particular, this implies that the fibers $M_C$ and $\tildeM_C$ are isomorphic.

In \cite[Proposition 2.1.4]{R1}, Rapagnetta shows that the locus $M_C^{lf} \subset M_C$ parameterizing sheaves which are pushforwards of line bundles on $C$ is open, dense, and non-canonically isomorphic to $\ker(a)$. Rapagnetta's proof also immediately implies the analogous statement for $N_C$. Rapagnetta then studies $\ker(a)$ by noting that there is a commutative diagram
\begin{equation}\label{a, tildea commutative diagram}
\begin{tikzcd}
 \Pic^0(C) \arrow[r,"\nu^*"] \arrow[d,"a"] &\Pic^0(\tildeC) \arrow[dl, "\tildea"]\\
 J
\end{tikzcd}
\end{equation} 
where $\tildea \colon \Pic^0(\tildeC) \to J$ be the morphism on Albanese varieties induced by the composition $\tildei = i \circ \nu$. However, Rapagnetta does not explicitly compute the number of connected components of $\ker(a)$ when $C \in N(2)$. We now give a slight generalization of Rapagnetta's argument, based on the argument in \cite[Remark 5.1 (5)]{MRS1}, which allows us to make this computation. 

Recall that the linear system $|2\theta|$ induces a double cover $\phi \colon J \to K$ onto a Kummer quartic $K \subset |2\theta|^\vee$ ramified at the 2-torsion points $J[2]$. Restricting $\phi$ to $C$, we see that $C$ is a double cover of its image $D = \phi(C)$ ramified at the nodes of $C$. There is a double cover $\tildephi \colon \tildeC \to \tildeD$ between normalizations ramified at the $2k$-points lying over the nodes of $C$. There are natural covering involutions $\iota$ on $C$ and $\tildeiota$ on $\tildeC$ which give rise to involutions $\sigma \coloneqq \iota^*$ on $\Pic^0(C)$ and $\tildesigma \coloneqq \tildeiota^*$ on $\Pic^0(\tildeC)$.

\begin{lem}\label{ker(a) subset Fix(sigma)}
	With the notation above, we have $\ker(a) \subset \Fix(\sigma)$. 
	\begin{proof}
		By \cite[Remark 5.1(3)]{MRS1}, the $-1$ involution on $J$ induces an involution via pullback on the moduli space $\bfM_{(0,2\theta,-4)}$ whose fixed locus contains $M_{(0,2\theta,-4)}$. The claim then follows by considering the restriction of this involution to the fiber over $C$.
	\end{proof}
\end{lem}

The fixed locus of $\sigma$ can be understood using the theory of Prym varieties for the branched double covers developed by Mumford in \cite{M1} in the smooth case and extended to the nodal case by Beauville in \cite{B3}.

\begin{lem}\label{description of fix(sigma)}
	$\Fix(\sigma) = \ker(1-\sigma)$ is connected if $k = 0,1$ and has $2^{k-1}$ connected components if $k=2, 3$. Each connected component is isomorphic to a $(\bbC^*)^k$-bundle over $\tildephi^*\Pic^0(\tildeD)$.
		\begin{proof}
		Consider the commutative diagram with exact rows
		\begin{equation*}
		\begin{tikzcd}
		1 \arrow[r] & (\bbC^*)^k \arrow[r] \arrow[d] & \Pic^0(C) \arrow[r,"\nu^*"] \arrow[d, twoheadrightarrow, "1- \sigma"] &\Pic^0(\tildeC) \arrow[r] \arrow[d,twoheadrightarrow, "1-\widetilde{\sigma}"] &1\\
		1 \arrow[r] & \ker(\alpha) \arrow[r] & (1-\sigma)\Pic^0(C) \arrow[r,"\alpha"] & (1-\widetilde{\sigma})\Pic^0(\tildeC) \arrow[r] & 1
		\end{tikzcd}
		\end{equation*}
		where $\alpha$ is the restriction of $\nu^*$ to the subgroup $(1-\sigma)\Pic^0(C)$. We note that $(1-\widetilde{\sigma})\Pic^0(\tildeC)$ and $(1-\sigma)\Pic^0(C)$ are the Prym varieties associated to the branched covers $\tildephi$ and $\phi$ respectively and are both Abelian surfaces.  
		
		When $k=0$, $C = \tildeC$ and the double cover $\phi \colon C \to D$ is unramified. In this case, the pullback $\phi^* \colon \Pic^0(D) \to \Pic^0(C)$ is not injective and $\ker(\phi^*) \simeq \bbZ/2\bbZ$. Mumford proves that $\ker(1-\sigma) \simeq \phi^*\Pic^0(D)$	and it follows that $\ker(1-\sigma)$ is connected (see \cite[\S 2 Data II (vi)]{M1}). When $k = 1, 2,3$, the double cover $\tildephi \colon \tildeC \to \tildeD$ is ramified and the pullback $\tildephi^* \colon \Pic^0(\tildeD) \to \Pic^0(\tildeC)$ is injective. Mumford proves that $\ker(1-\widetilde{\sigma}) \simeq \tildephi^*\Pic^0(\widetilde{D}) \times (\bbZ/2\bbZ)^{2k-2}$ (loc. cit.). 
		
		Beauville proves that $\ker(\alpha)$ is a finite group of order $2^{k-1}$ (see \cite[Remark 3.6]{B3}) which implies that there is no nontrivial group homomorphism $\bbC^* \to \ker(\alpha)$. The snake lemma then implies that there is a short exact sequence
		\begin{equation*}
		1 \to (\bbC^*)^k \to \ker(1-\sigma) \to \tildephi^*\Pic^0(\widetilde{D}) \times (\bbZ/2\bbZ)^{k-1} \to 1,
		\end{equation*}
		which proves the lemma. 
	\end{proof}
\end{lem}

The following observation, which was communicated to us by Rapagnetta, will help us describe $\ker(a)$ when $C \in N(2)$. 

\begin{lem}\label{O(p1-p2) not in identity component cor}
	Consider $C \in N(2)$ with nodes at $p,q \in J[2]$. Denote by $p_1,p_2$ the points in $\tildeC$ lying over $p$. Then the degree 0 line bundle $\calO_{\tildeC}(p_1-p_2)$ is not pullback of a line bundle from $\tildeD$. 
	
	\begin{proof}
		The involution $\tildeiota$ defines a $\bbZ/2\bbZ$ action on $\tildeC$ whose 	quotient is $\tildeD$. Since $p_1, p_2$ are in the fixed locus of the involution $\tildeiota$ on $\tildeC$, the line bundle $\calO_{\tildeC}(p_1-p_2)$ is $\bbZ/2\bbZ$ invariant and hence, is $\bbZ/2\bbZ$-linearizable (cf. \cite[Remark 7.2]{Dol}). Moreover, there are only two possible $\bbZ/2\bbZ$ linearizations. The Kempf Descent Lemma (see \cite[Proposition 4.2.15]{HL}) says that $\calO_{\tildeC}(p_1-p_2)$ descends to a line bundle on $\tildeD$ if and only if there exists a $\bbZ/2\bbZ$-linearization of $\calO_{\tildeC}(p_1-p_2)$ such that the linearized action on the fibers of $\calO_{\tildeC}(p_1-p_2)$ over the fixed points $p_1,p_2,q_1,q_2$ of $\tildeiota$ is trivial. Considering $\calO_{\tildeC}(p_1-p_2)$ as the subsheaf of the sheaf of rational functions consisting of rational functions on $\tildeC$ vanishing at $p_2$ and having at most a simple pole at $p_1$, one can use local coordinates around $p_1,p_2,q_1,q_2$ in which the $\bbZ/2\bbZ$-action on $\tildeC$ is linear to check that there is no $\bbZ/2\bbZ$-linearization satisfying the conditions of Kempf's Descent Lemma.
	\end{proof}
\end{lem}

\begin{lem}\label{kera = fix sigma lemma}
	For $k=0,1,2,3$, $\ker(a) = \Fix(\sigma)$.
	\begin{proof}
		For $k=0,1,2,3$, the inclusion $\ker(a) \subseteq \Fix(\sigma)$ was discussed in Lemma \ref{ker(a) subset Fix(sigma)}. We will show the reverse inclusion case by case. If $k = 0,1$, then $\Fix(\sigma)$ is connected by Lemma \ref{description of fix(sigma)} which implies that $\ker(a) = \Fix(\sigma)$. If $k = 2$, noting that the identity component $\ker(\tildea)^0$ can be identified with $\tildephi^*(\Pic^0(\tildeD))$, Lemma \ref{O(p1-p2) not in identity component cor} implies that $\ker(\tildea)$, and hence $\ker(a)$, is not connected. Since $\Fix(\sigma)$ has exactly two connected components by Lemma \ref{description of fix(sigma)}, we conclude that $\ker(a) = \Fix(\sigma)$.	If $k = 3$, then Rapagnetta shows in \cite[Proposition 2.1.4]{R1} that $\ker(a)$ has four connected components. Since $\Fix(\sigma)$ also has four connected components by Lemma \ref{description of fix(sigma)}, we conclude that $\ker(a) = \Fix(\sigma)$. 
	\end{proof}
\end{lem}

\begin{prop}\label{irred components over integral curves prop} 
	If $C \in S \cup N(1) \cup N(2) \cup N(3)$ has $k$ nodes, then the fiber 
	\begin{enumerate}
		\item $G_C$ is isomorphic to $\tildephi^*\Pic^0(\widetilde{D})$ if $k= 0,1$ and is isomorphic to $\tildephi^*\Pic^0(\widetilde{D}) \times (\bbZ/2\bbZ)^{k-1}$ if $k=2,3$;
		
		\item $\tildeM_C$ is irreducible of dimension three if $k=0,1$, has $2^{k-1}$ irreducible components of dimension three if $k=2,3$, is reduced, Cohen-Macaulay, and has an open dense subset parameterizing line bundles of degree $2$;
		
		\item $N_C$ is irreducible of dimension three if $k=0,1$, has $2^{k-1}$ irreducible components of dimension three if $k=2,3$, is reduced, Cohen-Macaulay, and has an open dense subset parameterizing line bundles of degree $1$;
	\end{enumerate}	
	\begin{proof}
		The statements about $G_C$ and the irreducible components follow from the lemmas in this section. We now argue that $\tildeM_{C}$ is reduced. First recall that $\tildeM_C \simeq M_C$ since every point in $M_C$ is a stable sheaf. Since the Le Potier support morphism is smooth at every point corresponding to a line bundle by \cite[Proposition 2.8]{LP}, the fiber $M_{C}$ is reduced on the open dense subset $M_C^{lf}$. By Remark \ref{flatness of the Lagrangian fibrations}, the fiber $M_{C}$ is Cohen-Macaulay and thus, is reduced everywhere. The same argument shows that $N_C$ is reduced.
	\end{proof}
\end{prop}

\subsection{Monodromy of Irreducible Components over \texorpdfstring{$N(2)$}{N(2)}}\label{monodromy of irred comp section}
In this section, we study the monodromy of the irreducible components of the fibers of $\tildem$ and $n$ over $N(2)$. Recall that this is equivalent to studying the monodromy of the irreducible components of the fibers of $m$ and $n$ over $N(2)$. The strategy used in this section was communicated to us by Rapagnetta. 

Recall from Proposition \ref{stratification of |2theta|} that $\overline{N(2)} = \cup_{p,q \in J[2], p \neq q} \Npq$ where $\Npq \simeq \bbP^1$. Let $\Upq = \Npq \setminus NR \simeq \bbC^*$, and $\Vpq = \Npq \setminus (NR \cup N(3)) \simeq \bbP^1 \setminus \{8\text{ points}\}$. In particular, $N(2) = \sqcup_{p,q \in J[2], p \neq q} \Vpq$. In the remainder of the section, we work with a fixed choice of $p$ and $q$. Consider the commutative diagram
\begin{equation}\label{kummer k3 diagram}
\begin{tikzcd}
\tildeJ \arrow[r, "\tildephi"] \arrow[d,"\varepsilon"] & \tildeK \arrow[d, "\tau"]\\
J \arrow[r, "\phi"] & K 
\end{tikzcd}
\end{equation}
where $\tildeJ$ is the blow up of $J$ along the locus of 2-torsion points $J[2]$, $K \simeq J / \pm 1$, and $\tildeK$ is the Kummer K3 surface associated to the Abelian surface $J$. Denote the exceptional divisors of $\varepsilon$ and $\tau$ by $F_r$ and $E_r$ respectively for $r \in J[2]$.

Let $\calC_{pq} \subset J \times \Vpq$ be the family of integral curves in $\Vpq$, $\tildecalC_{pq} \subset \tildeJ \times \Vpq$ be the family of strict transforms, and $\tildecalD_{pq} \coloneqq (\tildephi \times id)(\tildecalC_{pq}) \subset \tildeK \times \Vpq$ be the corresponding image. Using Equation \ref{kummer k3 diagram}, one can see that if $C_b$ is the fiber of $\calC_{pq}$ over $b \in \Vpq$, then the fiber of $\tildecalC_{pq}$ over $b$ is the normalization $\tildeC_b$ of $C_b$ and the fiber of $\tildecalD_{pq}$ over $b$ is the normalization $\tildeD_b$ of $D_b \coloneqq \phi(C_b)$. There are natural 2-sections $S_p \coloneqq (F_p \times \Vpq) \cap \tildecalC_{pq} \subset \tildecalC_{pq}$ and $T_p \coloneqq (E_p \times \Vpq) \cap \tildecalD_{pq} \subset \tildecalD_{pq}$. Some geometry of these 2-sections is discussed in the following two lemmas.

\begin{lem}\label{(x,b) in Tp lemma}
	Given $(x,b) \in T_p$, there exists a unique line $\ell_x \subset |2\theta|^\vee$ and a unique plane $H_b \subset |2\theta|^\vee$ such that $\ell_x \subset H_b$ and $\phi(q) \notin \ell_x$. 
	\begin{proof}
		For $b \in \Vpq$, let $C_b$ be the corresponding curve in $J$ with nodes at $p,q \in J[2]$. Then $D_b = \phi(C_b) = H_b \cap K$ for a unique plane $H_b \subset |2\theta|^\vee$ containing $\phi(p)$ and $\phi(q)$.
		Now given $(x,b) \in T_p$, we see that $x \in \tildeD_b \cap E_p$. In particular, $x$ corresponds to a tangent direction to $D_b$ at the node $\phi(p)$. Let $\ell_x$ be the unique line in $|2\theta|^\vee$ corresponding to this tangent direction. Since $D_b = H_b \cap K$, we see that $\ell_x \subset H_b$.
		
		To see that $\phi(q) \notin \ell_x$, note that $D_b \subset H_b$ is a quartic curve which is nodal at $\phi(p)$ and $\phi(q)$. Since $\ell_x$ is tangent to $D_b$ at the node $\phi(p)$, the intersection multiplicity of $\ell_x \cap D_b$ at $\phi(p)$ is 3. If $\phi(q) \in \ell_x$, then the intersection multiplicity of $\ell_x \cap D_b$ at $\phi(p)$ would be strictly larger than 1 since $D_b$ is singular at $\phi(p)$. This would imply that the intersection number $\ell_x \cdot D_b$ would be larger than 5 which contradicts the fact that $D_b$ is a quartic curve in $H_b$.
	\end{proof}
\end{lem}

\begin{lem}
	The 2-section $S_p$ is connected. 
	\begin{proof}
		Since the map $\tildephi \colon \tildeJ \to \tildeK$ in Equation \ref{kummer k3 diagram} is an isomorphism when restricted to the exceptional divisor $F_p$, $S_p \simeq T_p$. Since $S_p \simeq T_p$ has no zero dimensional components, it suffices to show that the projection map $T_p \to E_p$ is injective as then $S_p \simeq T_p$ would be isomorphic to a Zariski open subset of the irreducible conic $E_p$. To see the desired injectivity, fix a hyperplane $H$ containing $\phi(q)$ but not $\phi(p)$. Given any $(x,b) \in T_p$, let $\ell_x$, $H_b$ be as in Lemma \ref{(x,b) in Tp lemma}. Let $\ell_b \subset H$ be the line $H_b \cap H$ and let $y$ be the point $\ell_x \cap H$. Since $\phi(q) \in H_b$ for all $b$ and $\phi(q) \in H$ by choice of $H$, $\phi(q) \in \ell_b$. Since $y$ and $\phi(q)$ are distinct by Lemma \ref{(x,b) in Tp lemma}, $\ell_b$ is the unique line through $y$ and $\phi(q)$. Now if $(x,b), (x',b') \in T_p$ with $x = x'$, then $y=y'$. It follows that the lines $\ell_b$ and $\ell_b'$ coincide which implies that $b=b'$.
	\end{proof}
\end{lem}

For any integer $d$, let $\Pic^d_{\calC_{pq}/\Vpq}$ and $\Pic^d_{\tildecalC_{pq}/\Vpq}$ denote the respective relative Picard varieties. In view of Section \ref{irred components over integral curves section}, the loci $N_{\Vpq}^{lf}\subset N_{\Vpq}$ and $M_{\Vpq}^{lf} \subset M_{\Vpq}$ parameterizing sheaves which are the pushforwards of line bundles are dense and can be identified with the fiber over $0$ of the map $\bfa \colon \Pic^d_{\calC_{pq}/\Vpq} \to J$, which fiberwise is the map defined in Equation \ref{a equation}, for $d = 1,2$ respectively. To study the monodromy of the components of $N$ and $M$ over $\Vpq$, we will study the connectivity of $\bfa\inv(0)$ for $d = 1,2$. In view of Equation \ref{a, tildea commutative diagram}, the diagram
\begin{equation*}\label{albanese diagram over Vpq}
\begin{tikzcd}
\Pic^d_{\calC_{pq}/\Vpq} \arrow[d, twoheadrightarrow, "\bfa"] \arrow[r,twoheadrightarrow] &\Pic^d_{\tildecalC_{pq}/\Vpq} \arrow[dl, twoheadrightarrow, "\tildebfa"]\\
J 
\end{tikzcd}
\end{equation*}
commutes. In particular, $\bfa \inv(0)$ is connected if and only if $\tildebfa \inv(0)$ is connected.

\begin{prop}\label{M,N over Vpq prop}
	$N_{\Vpq}$ is irreducible and $M_{\Vpq}$ is reducible. 
	\begin{proof}
		We first discuss $N_{\Vpq}$. By the discussion above, it suffices to show that the fiber of the map $\tildebfa \colon \Pic^1_{\tildecalC_{pq}/\Vpq} \to J$ over $0$ is connected. Consider the special case where $p = 0 \in J[2]$. Let $b \in \Vpq$ and $\tildeC_b$ be the fiber of the family of curves $\tildecalC_{pq}$ at $b$. Let $\tildebfa_b \colon \Pic^1(\tildeC_b) \to J$ be the restriction of $\tildebfa$ and recall the fiber of $\tildebfa_b$ over $0 \in J$ consists of two connected components by Proposition \ref{irred components over integral curves prop}. Let $p_1,p_2$ be the two points in $\tildeC_b$ lying over the node $p \in C_b$. Since $p_1, p_2$ lie over $p = 0$, Grothendieck-Riemann-Roch for the composition $\tildeC_b \to C_b \to J$ implies that $\tildebfa_b(\calO_{\tildeC_b}(p_i)) = 0$ for $i = 1,2$. 
		
		Recall that by Corollary \ref{O(p1-p2) not in identity component cor}, the degree 0 line bundle $\calO_{\tildeC_b}(p_1-p_2)$ is not in the identity component of kernel of the map $a \colon \Pic^0(\tildeC_b) \to J$. This implies that the degree 1 line bundles $\calO_{\tildeC_b}(p_1)$ and $\calO_{\tildeC_b}(p_2)$ must lie in different connected components of $\tildebfa_b \inv(0)$. Consider the 2-section $S_p \subset \tildecalC_{pq}$ and let $\widehat{\calC}_{pq} \coloneqq \tildecalC_{pq} \times_{\Vpq} S_p$ be the fiber product. The preimage of the 2-section $S_p \subset \tildecalC_{pq}$ in $\widehat{\calC}_{pq}$ splits into two irreducible components $S_1$ and $S_2$.	Using representability of the Picard functor, the line bundle $\calO_{\widehat{\calC}_{pq}}(S_1)$ on $\calC_{S_p}$ defines a section $s \colon S_p \to \Pic^1_{\tildecalC_{pq}/\Vpq}$ such that $s(p_i,b) = \calO_{\tildeC_b}(p_i)$. Since $s(p_1,b)$ and $s(p_2,b)$ are in different connected components of $\tildebfa \inv(0)$ and $S_p$ is connected, $\tildebfa \inv(0)$ is connected. 
		
		If $p \neq 0 \in J[2]$, consider the translation by $p$ morphism $t_p \colon J \to J$. Pullback by $t_p$ determines an automorphism of the linear system $|2\theta|$ sending the line $\Npq$ to the line $N_{0r}$ where $r=p+q$. The connectivity of $\tildebfa \inv(0)$ then follows from the $p = 0$ case. 
		
		We now discuss $M_{\Vpq}$. By the discussion above, it suffices to show that the fiber of the map  $\tildebfa \colon \Pic^2_{\tildecalC_{pq}/\Vpq} \to J$ over $0$ is disconnected. As in the previous case, let $\tildebfa_b \colon \Pic^2(\tildeC_b) \to J$ be the restriction of $\tildebfa$ for $b \in \Vpq$. Since the fiber of $\tildebfa_b$ over $0 \in J$ consists of two connected components by Proposition \ref{irred components over integral curves prop}, it suffices to exhibit a section of $\Pic^2_{\tildecalC_{pq}/\Vpq} \to \Vpq$ which factors through $\bfa\inv(0)$. To see that such a section exists, consider the projection $\widetilde{\pi}_{pq} \colon \tildecalC_{pq} \to \tildeJ$ and the line bundle $\widetilde{\pi}_{pq}^*\calO_{\Jhat}(F_p)$ on $\tildecalC_{pq}$. Using representability of the Picard functor, this line bundle defines a section $s \colon \Vpq \to \Pic^2_{\tildecalC_{pq}/\Vpq}$ such that $s(p_i,b) = \calO_{\tildeC_b}(p_1+p_2)$. Since $p_1, p_2$ lie over the 2-torsion point $p \in J[2]$, we see that this section factors through $\tildebfa \inv(0)$.
	\end{proof} 
\end{prop}

\subsection{Irreducible Components Over Non-Integral Curves}\label{irred components over non-integral curves section}
We begin with the following proposition on the fibers of $\tildem$ and $n$ over the locus of reduced, but reducible curves, i.e. over $R(1) \cup R(2)$.

\begin{prop}\label{irred comp of tildeM, N over R(1) u R(2)}
	If $C \in R(1) \cup R(2)$, then the fiber
	\begin{enumerate}
		\item $G_C$ is isomorphic to a $\bbC^*$-bundle over the Abelian surface $J$ if $C \in R(1)$ and is isomorphic to a $\bbC$-bundle over $J$ if $C \in R(2)$;
		\item $\tildeM_C$ has two irreducible components of dimension three and has a non-dense open subset parameterizing line bundles of bidegree $(1,1)$;
		\item $N_C$ has two irreducible components of dimension three, is reduced, and has an open dense subset parameterizing line bundles of bidegree $(0,1)$ and $(1,0)$.
	\end{enumerate}
	\begin{proof}
	The first claim about the group scheme is proven by Rapagnetta in \cite[Proposition 2.1.4]{R1}. To see the second claim, recall that $\tildem$ is the composition of the symplectic resolution $\pi \colon \tildeM \to M$ and the support morphism $m \colon M \to B$. Rapagnetta's proof in \cite[Proposition 2.1.4]{R1} also shows that the fiber $M_C$ is irreducible with an open dense subset parameterizing line bundles with bidegree $(1,1)$ and that the strictly semi-stable locus in $M_C$ is isomorphic to the Abelian surface $J$. By Remark \ref{fibers of symplectic resolution},  $\tildeM_C$ consists of two irreducible components, namely the strict transform of $M_C$ under the symplectic resolution and the exceptional divisor restricted to the fiber over $C$, which is isomorphic to a $\bbP^1$-bundle over $J$. A straightforward adaptation of Rapagnetta's argument in \cite[Proposition 2.1.4]{R1} gives the third claim. The claim about reducedness of $N$ follows from the same argument as in the proof of Proposition \ref{irred components over integral curves prop}.
	\end{proof}
\end{prop}

We now study the irreducible components of the fibers of the Lagrangian fibrations $\tildem$ and $n$ over the locus of non-reduced curves, i.e. over the stratum $NR$. We remark that the analysis in this section is simpler than the analysis over the non-reduced locus in the $OG10$ case since the underlying reduced curves in our case are always smooth.  

The sheaves parameterized by the fibers $M_C$ and $N_C$ can be divided into two types, namely sheaves of type I and type II (cf. \cite[\S 3.7.2]{dCRS}). A sheaf $\scrF$ on $J$ with Fitting support $C = 2\Cred$ is said to be of type $I$ if the composition of the natural morphisms $\calO_J \to \calO_C \to \sheafEnd_J(\scrF)$	factors via the natural surjection $\calO_J \to \calO_{\Cred}$ and is said to be of type $II$ otherwise. In particular, sheaves of type $I$ are sheaves on $J$ with schematic support $\Cred$ and sheaves of type $II$ are sheaves on $J$ with schematic support $C$. 

\begin{lem}\label{type I sheaf}
	If $C = 2\Cred \in NR$, then
	\begin{enumerate}
		\item $M_C^I$, with its reduced induced structure, is irreducible and isomorphic to $\bbP^3$. Moreover, the locus of strictly semi-stable bundles in $M_{C}^{I}$, with its reduced induced structure, is isomorphic to a Kummer quartic surface,
		
		\item $N_C^I$, with its reduced induced structure, is irreducible and is the intersection of two quadrics in $\bbP^5$. 
	\end{enumerate}
	\begin{proof}
		The sheaves parameterized by $M_C^I$ (resp. $N_C^I$) are pushforwards of  semi-stable rank two vector bundles on $\Cred$ of degree $0$ (resp. degree $-1$) and fixed determinant. The structure of such moduli spaces is described in \cite{NS}.
	\end{proof}
\end{lem}

\begin{lem}\label{type II sheaf}
	If $C = 2\Cred \in NR$, then 
	\begin{enumerate}
		\item $M_{C}^{II}$ is reduced and consists of 16 connected components, each of which is isomorphic to $\bbC^3$,
		\item $N_C^{II}$ is connected.
	\end{enumerate}
	\begin{proof}
		The first statement was proven in \cite[Proposition 2.1.4]{R1}. Adapting Rapagnetta's argument, one can show that sheaves $\scrF$ parameterized by $N_C^{II}$ are precisely the extensions
		\begin{equation*}
		0 \to i_{\red*}(F \otimes \calI \otimes \calO_{\Cred}(x)) \to \scrF \to i_{\red*}(F) \to 0
		\end{equation*}
		where $F = \scrF|_{\Cred}/Tors$ is a degree $0$ line bundle on $\Cred$, $x \in \Cred$ is a point, and $\sum c_2(\scrF) = 0$. The condition $\sum c_2(\scrF) = 0$ implies that $0 = \sum c_2(i_{\red*}(F^{\otimes 2} \otimes \calI)) - i_{\red}(x)$. The collection of line bundles $F$ on $\Cred$ which satisfy this equation is precisely the fiber over $0 \in J$ of the map	$g \colon \Cred \times \Pic^{0}(\Cred) \to J$ sending the pair $(x,F)$ to the point $\sum c_2(i_{\red*}(F^{\otimes2} \otimes \calI) - i_{\red}(x) \in J$. The fiber $g \inv (0)$ can be described by the fiber product diagram
		\begin{center}
			\begin{tikzcd}
				g\inv(0) \arrow[r] \arrow[d] \arrow[dr,phantom,"\square"] &\Pic^{0}(\Cred) \arrow[d, "h"] \\
				\Cred \arrow[r,"i_{\red}"] &J
			\end{tikzcd}
		\end{center}
		where the map $h$ is a degree 16 isogeny sending a degree 0 line bundle $F$ to the point $\sum c_2(i_{\red*}(F^{\otimes2} \otimes \calI)) \in J$ and $i_{\red} \colon \Cred \to J$ can be identified with the Albanese morphism. We conclude that $g\inv(0)$ is connected since the pullback of a connected \'etale cover of the Albanese variety by the Albanese map is connected (cf. \cite[Remark V.14 (5)]{B1}). For a fixed line bundle $F$ in $g \inv(0)$, the argument in \cite[Proposition 3.7.19]{dCRS} implies that the possible extensions are parameterized by $\bbP(\Ext^1_C(F,K)) \setminus \bbP(\Ext^1_{\Cred}(F,K)) \simeq \bbC^2$. It follows that $N_C^{II}$, with its reduced, induced structure, is isomorphic to a $\bbC^2$-bundle over $g\inv(0)$. 
	\end{proof}
\end{lem}

\begin{prop}\label{irred comp of tildeM, N over NR}
	If $C = 2\Cred \in NR$, then the fiber
	\begin{enumerate}
		\item $G_C$ has 16 connected components, each of which is isomorphic to $\bbC^3$;
		\item $\tildeM_C$ consists of 34 irreducible components of dimension three and has an open set, which is not dense, parameterizing pushforwards of degree $1$ line bundles from $C$;
		\item $N_C$ has two irreducible components of dimension three. Moreover, no sheaf in $N_C$ is the pushforward of a line bundle on $C$.
	\end{enumerate}
	\begin{proof}
		The first claim about the group scheme is proven by Rapagnetta in \cite[Proposition 2.1.4]{R1}. To see the second claim, recall that by Lemma \ref{type I sheaf}, the strictly semi-stable locus $M_{C}^{ss} = M_{C} \cap \Sigma$, with its reduced induced structure, is isomorphic to a Kummer quartic surface $K \subset |2\theta|$. Under this isomorphism, $M_C \cap \Omega$ consists of the 16 nodes of $K$. It follows from Remark \ref{fibers of symplectic resolution} that $\tildeM_{C}$ consists of $34$ irreducible components, namely the strict transforms of $M_C^I$ and the 16 components in $M_{C}^{II}$ and the 17 irreducible components which constitute the exceptional divisor of the resolution restricted to the fiber over $C$. The third claim about $N_C$ follows from the two lemmas above.
	\end{proof}
\end{prop}

\section{The Top Degree Direct Image Sheaves}\label{top degree direct image sheaves section}

\begin{fact}\label{trace lemma rmk}
	Let $f \colon X \to T$ be a flat morphism of relative dimension $d$. Then
	\begin{enumerate}
		\item\label{fact 1} there is a trace morphism $\Tr_f \colon R^{2d} f_! \bbQ_X(d) \to \bbQ_T$ which is an isomorphism if and only if all of the fibers of $f$ have a unique irreducible component of dimension $d$ (cf. \cite[Th\'eor\`eme 2.9 and Remarque 2.10.1]{SGA4}),
		\item\label{fact 2} if $f$ has reduced fibers, then the sheaf $R^{2d}f_!\bbQ_X(d)$ is the $\bbQ$-linearization of the sheaf of sets of irreducible components of the fibers of $f$ (cf. \cite[Lemme 7.1.8]{N1}). 
	\end{enumerate}
\end{fact}

Recall from Section \ref{monodromy of irred comp section} that $N(2) = \sqcup_{p,q \in J[2], p \neq q} \Vpq$ where $\Vpq = \Npq \setminus (NR \cup N(3)) \simeq \bbP^1 \setminus \{8\text{ points}\}$.  Let $R^0 \coloneqq R(1) \cup R(2) = R \setminus NR$.

\begin{prop} \label{direct image sheaves in top degree prop}
	Define $R^6_{\tildeM} \coloneqq R^6\widetilde{m}_*\bbQ_{\tildeM}$ and $R^6_N \coloneqq R^6n_*\bbQ_N$. There are isomorphisms of constructible sheaves: 
	\begin{equation}\label{highest direct image sheaf over S u N(1)}
		R^6_{\tildeM}|_{S \cup N(1)} \simeq R^6_N|_{S \cup N(1)} \simeq \bbQ_{S \cup N(1)}
	\end{equation}
	
	\begin{equation}\label{highest direct image sheaf over N(3) and NR}
		R^6_{\tildeM}|_{N(3)} \simeq R^6_N|_{N(3)} \simeq \bbQ_{N(3)}^{\oplus 4}; \quad R^6_{\tildeM}|_{NR} \simeq\bbQ_{NR}^{\oplus 34}; \quad R^6_N|_{NR} \simeq \bbQ_{NR}^{\oplus 2}
	\end{equation}
	
	\begin{equation}\label{highest direct image sheaf of tildeM, N over R}
		R^6_{\tildeM}|_{R^0} \simeq \bbQ_{R^0}^{\oplus 2}; \quad R^6_N|_{R^0}\simeq \bbQ_{R^0} \oplus \scrL_{R^0}
	\end{equation}
	
	\begin{equation}\label{highest direct image sheaf of tildeM, N over Vpq}
		R^6_M|_{\Vpq} \simeq R^6_{\tildeM}|_{\Vpq} \simeq \bbQ_{\Vpq}^{\oplus 2}; \quad R^6_N|_{\Vpq} \simeq \bbQ_{\Vpq} \oplus \scrL_{\Vpq},
	\end{equation}
		where $\scrL_{R^0}$ is the rank one local system on $R^0$ corresponding to the \'etale double cover $J \setminus J[2] \to R^0$ sending a point $x$ to the curve $\theta_x + \theta_{-x}$ and $\scrL_{\Vpq}$ is a rank one local system on $\Vpq$ corresponding to $\pi_1(\Vpq)$ representation for which loops around $N(3)$-points act by the identity and loops around $NR$-points act by $-1$.
	\begin{proof}
	Recall that by Remark \ref{flatness of the Lagrangian fibrations}, the morphisms $m,\tildem$ and $n$ are all flat. We begin by proving Equation \ref{highest direct image sheaf over S u N(1)}. Since the fibers of $M, \tildeM$, and $N$ over the loci $S \cup N(1) \subset B$ are integral by Proposition \ref{tildeM,M,N irred components prop}, Fact \ref{fact 1} implies the equality. The isomorphisms in Equation \ref{highest direct image sheaf over N(3) and NR} follow immediately from Proposition \ref{tildeM,M,N irred components prop} since $N(3)$ and $NR$ are zero dimensional.
	
	We next prove the left isomorphism in Equation \ref{highest direct image sheaf of tildeM, N over R}. Proposition \ref{DT for symplectic resolution} on the Decomposition Theorem for the symplectic resolution $\pi \colon \tildeM \to M$ and proper base change imply that $R\tildem_* \bbQ_{\tildeM_{R^0}} \simeq Rm_* \bbQ_{M_{R^0}} \oplus Rr_* \bbQ_{\Sigma_{R^0}}[-2](-1)$, where $r$ denotes the restriction of the support morphism $m \colon M \to B$ to $\Sigma_{R^0}$. Thus, $R^6\tildem_* \bbQ_{\tildeM_{R^0}} \simeq R^6m_* \bbQ_{M_{R^0}} \oplus R^4r_* \bbQ_{\Sigma_{R^0}}(-1)$. Since the fiber $M_C$ is irreducible for any $C \in R^0$ by Proposition \ref{irred comp of tildeM, N over R(1) u R(2)}, Fact \ref{fact 1} implies that $R^6m_* \bbQ_{M_{R^0}} \simeq \bbQ_{R^0}$. By Lemma \ref{type I sheaf}, the fibers of the map $\Sigma_{R^0} \to R^0$ are isomorphic to the Abelian surface $J$ and hence irreducible. Fact \ref{fact 1} then implies that $R^4r_* \bbQ_{\Sigma_{R^0}} \simeq \bbQ_{R^0}$ and the desired isomorphism follows. 
	
	We next prove the right isomorphism in Equation \ref{highest direct image sheaf of tildeM, N over R}. Recall that by Proposition \ref{irred comp of tildeM, N over R(1) u R(2)}, the morphism $N_{R^0} \to R^0$ is flat, with reduced fibers each having two irreducible components. Fact \ref{fact 2} implies that $R^6\bbQ_{N_{R^0}}$ is the $\bbQ$-linearization of the sheaf of sets $\Irr(N_{R^0})$ of irreducible components of $N_{R^0}$ which is locally constant with stalks of cardinality two. Let $N_{R^0}^{lf} \to R^0$ be the locus parameterizing sheaves which are the pushforwards of line bundles. Again by Proposition \ref{irred comp of tildeM, N over R(1) u R(2)}, the morphism $N_{R^0}^{lf} \to R^0$ is smooth, dense in every fiber, and surjective with two connected components which parameterize line bundles of bi-degree $(0,1)$ and bi-degree $(1,0)$. It follows that the sheaf of sets $\Irr(N_{R^0})$ can be identified with the sheaf of connected components of $N_{R^0}^{lf}$. Using this interpretation, we now examine the monodromy of $\Irr(N_{R^0})$. The map $J \setminus J[2] \to R^0$ sending a point $x$ to $\theta_x + \theta_{-x}$ is an \'etale double cover which induces a monodromy action of $\pi_1(R^0)$ on $J \setminus J[2]$ interchanging the fibers of the covering map. It follows that the monodromy action on the family of curves over $R^0$ swaps the components of the broken curves $\theta_x + \theta_{-x}$. Thus, the bidegrees of the line bundles are swapped as well and it follows that the components of $N_{R^0, C}^{lf}$ are swapped. Since $(R^6n_*\bbQ_N)|_{R^0}$ is the $\bbQ$-linearization of the sheaf of sets $\Irr(N_{R^0})$, the desired isomorphism follows.

	We next prove the left isomorphism in Equation \ref{highest direct image sheaf of tildeM, N over Vpq}. Recall that by Proposition \ref{irred comp of tildeM, N over NR}, the morphism $M_{\Vpq} \to \Vpq$ is flat, with reduced fibers each having two irreducible components. Fact \ref{fact 2} implies that $R^6\bbQ_{M_{\Vpq}}$ is the $\bbQ$-linearization of the sheaf of sets $\Irr(M_{\Vpq})$ of irreducible components of $M_{\Vpq}$. Proposition \ref{M,N over Vpq prop} implies that the sheaf $\Irr(M_{\Vpq})$ is the constant sheaf with stalks of cardinality 2 and the desired isomorphism follows.	
	
	We finally prove the right isomorphism in Equation \ref{highest direct image sheaf of tildeM, N over Vpq}.  Recall that by Proposition \ref{irred comp of tildeM, N over NR}, the morphism $N_{\Vpq} \to \Vpq$ is flat, with reduced fibers each having two irreducible components. Fact \ref{fact 2} implies that $R^6\bbQ_{N_{\Vpq}}$ is the $\bbQ$-linearization of the sheaf of sets $\Irr(N_{\Vpq})$ of irreducible components of $N_{\Vpq}$. Proposition \ref{M,N over Vpq prop} implies that $\Irr(N_{\Vpq})$ must be a non-trivial locally constant sheaf with stalks of cardinality two. 
	
	We now claim that loops around $N(3)$ points do not interchange the two irreducible components of a fiber $N_C$. To see this, note that Proposition \ref{irred components over integral curves prop} implies that the morphism $N_{\Vpq} \to \Vpq$ can be extended to a morphism $N_{\Upq} \to \Upq$ which is also flat with reduced fibers. Moreover, the locus $N_{\Upq}^{lf}$ parameterizing sheaves which are pushforwards of line bundles is dense in every fiber. Now fix any $\scrF \in N_C^{lf}$ where $C \in N(3) \cap \Npq = \Upq \setminus \Vpq$. Since the Le Potier support morphism is smooth at every point of $N_C^{lf}$, the map $n$, viewed in the analytic topology, is a submersion at $\scrF$. It follows that there exists a small disk $\Delta_{pq} \subset \Upq$ about $C \in N(3)$ and a smooth local section $s \colon \Delta_{pq} \to N^{lf}_{\Delta_{pq}}$ such that $s(0) = \scrF$. Restricting this section to $\Delta_{pq}^*$ gives a smooth local section of $N^{lf}_{\Delta_{pq}^*}$ which implies that $N^{lf}_{\Delta{pq}^*}$ is disconnected. In particular, loops around $N(3)$ points do not interchange the two irreducible components of a fiber.
	
	Since $\Irr(M_{\Vpq})$ is a non-trivial locally constant sheaf and the loops around $N(3)$ points act trivially on the irreducible components, it extends to a non-trivial locally constant sheaf on $\Vpq = \Npq \setminus NR \simeq \bbC^*$. We conclude that loops around both $NR$ points must interchange the two irreducible components of a fiber and the desired isomorphism follows. 
\end{proof}
\end{prop}

\section{Ng\^o Strings}\label{Ngo strings section}
In this section, we will use our knowledge of the direct image sheaves in top degree for the Lagrangian fibrations $\tildem$ and $n$ and the Ng\^o Support Theorem to determine the Decomposition Theorems for $R\tildem_*\bbQ_{\tildeM}$ and $Rn_*\bbQ_N$ in Proposition \ref{DT for tildeM and N}. We first introduce several relevant strings that will appear in the Decomposition Theorems for $R\tildem_*\bbQ_{\tildeM}$ and $Rn_*\bbQ_N$.

\subsection{The Relevant String over \texorpdfstring{$B$}{B}}\label{strings over B section}
Consider the group scheme $g \colon G \to B$ described in Equation \ref{group scheme def}. Over the locus $S \subset B$ of smooth curves, the map $g \colon G_S \to S$ is smooth and the fibers of $G$ are 3-dimensional Abelian varieties by Proposition \ref{irred components over integral curves prop}. Consider the higher direct image sheaves $\Lambda^i_B \coloneqq R^ig_*\bbQ_{G_S}$. We introduce the following complex, viewable in $D^bMHM_{alg}(B)$ or $D^b(B,\bbQ)$, which we will call a string:
\begin{equation}
\scrI_B \coloneqq \oplus_{i=0}^6 \scrIC_B(\Lambda^i_B)[-i]. 
\end{equation}

\begin{lem}\label{cohomology sheaf of I_B}
	The cohomology sheaf $\calH^1(\scrIC_B(\Lambda^5_B))$ vanishes. In particular, $\calH^6(\scrI_B) = \bbQ_B \oplus \calH^2(\scrIC_B(\Lambda_B^4))$	where $\calH^2(\scrIC_B(\Lambda_B^4))$ is a skyscraper sheaf.
	\begin{proof}
		Let $\gamma \colon \calC \to B$ be the universal family of curves in the linear system $B$. Let $\calC_S$ be the restriction of the family over the locus of smooth curves and let $\Lambda^1_\gamma \coloneqq R^1\gamma_*\bbQ_{\calC_S}$. Analyzing the Decomposition Theorem for $\gamma$ shows that $\calH^1(\scrIC_B(\Lambda^1_{\gamma})) = 0$. The natural closed embedding $G_S \to \Pic^0_{\calC/S}$ induces a surjective morphism of semi-simple local systems $\Lambda^1_\gamma \to R^1g_*\bbQ_{G_S}$.  Semi-simplicity then implies that $\scrIC_B(\Lambda^1_B)$ is a direct summand of $\scrIC_B(\Lambda^1_{\gamma})$.  The claim then follows from Relative Hard Lefschetz. The cohomology sheaf $\calH^2(\scrIC_B(\Lambda_B^4))$ is a (possibly trivial) skyscraper sheaf by the support conditions for intersection complexes.
	\end{proof}
\end{lem}

\subsection{The Relevant Strings over \texorpdfstring{$R$}{R}}\label{strings over R section}
Throughout this section, we use the canonical identification $J^\vee \simeq \Pic^0(C_0)$ where $C_0$ is our fixed genus 2 curve. Let $A \coloneqq J^\vee \times J$. There is a commutative diagram
\begin{equation}\label{relevant fibrations over R}
\begin{tikzcd}
A \arrow[r, "pr_2"] \arrow[d,"q"] \arrow[dr, "c"] &J \arrow[d,"q"]\\
\Sigma \arrow[r,"r"] & R
\end{tikzcd}
\end{equation}	
where $pr_2 \colon A \to J$ is projection onto the second factor, the map $q \colon A \to \Sigma$ sends $(L,x)$ to $[i_{x*}L \oplus i_{-x*}L^\vee]$ where $i_{\pm x} \colon C_0 \to J$ is the embedding of $C_0$ into $J$ with image $\theta_{\pm x}$, the map $q \colon J \to R$ sends $x$ to $\theta_x + \theta_{-x}$, and $r$ is restriction of the support morphism $m \colon M \to B$ to $\Sigma$. Note that both $q$  maps are double covers ramified along the locus of 2-torsion points. We introduce the following complexes, viewable in either $D^bMHM_{alg}(R)$ or $D^b(R,\bbQ)$, which we will again call strings. Let $R^0 \coloneqq R(1) \cup R(2) = R \setminus NR$, $\Lambda_{R^0}^i \coloneqq R^ir_*\bbQ_{\Sigma_{R^0}}$,

\begin{equation}
\scrI^+_{R} \coloneqq \oplus_{i=0}^4 \scrIC_R(\Lambda_{R^0}^i)[-i],
\end{equation}
\begin{equation}
\scrI^-_{R} \coloneqq \oplus_{i=0}^4 \scrIC_R(\Lambda_{R^0}^i \otimes \scrL_{R^0})[-i],
\end{equation}
where $\scrL_{R^0}$ is the rank one local system on $R^0$ associated to the \'etale double cover $J_{R^0} \to R^0$. We will see in Section \ref{DT for tildeM, N section} that the Ng\^o Support Theorem implies that the local systems appearing in the string $\scrI_R^+$ can be identified with the cohomology of the Abelian part of the identity component of the group scheme $G$ over the Zariski open subset $R^0 \subset R$ .

\begin{lem}\label{DT for A/+-1 to R}
	There is an isomorphism $Rr_*\bbQ_{\Sigma} \simeq \scrI_R^+$. Moreover, the intersection complexes $\scrIC_R(\Lambda^1_{R^0})$ and $\scrIC_R(\Lambda^3_{R^0})$ are sheaves.
	\begin{proof}
		By looking at the regular part of $r$, $\scrI_R^+[4]$ is a direct summand of $Rr_*\bbQ_{\Sigma}[4]$. Using the fact that all fibers of $r$ are irreducible, one can show that there are no additional summands in perverse degrees $1$ and $2$. Relative Hard Lefschetz then implies that the only possible  additional summands are intersection complexes supported on one dimensional subvarieties of $R$ placed in perverse degree $0$. Such a summand would contribute non-trivially to $R^3r_*\bbQ_{\Sigma}$ in codimension one. However, since there are only finitely many singular fibers of $r$, no such summand can exist. We conclude that there are no additional summands appearing in the Decomposition Theorem for $r$. 
		
		The claim about the intersection complexes being sheaves follows immediately from $r$ having all irreducible fibers and Relative Hard Lefschetz. 
	\end{proof}
\end{lem}

We now turn our attention to the map $c = r \circ q \colon A \to R$ described in Equation \ref{relevant fibrations over R}. We note that the fibers of $c$ are straightforward to describe. Using the top right triangle in Equation \ref{relevant fibrations over R}, we see that fiber over a point in $R^0$ consists of two irreducible components, each of dimension two, while the fiber over a point in $NR$ is irreducible of dimension two. 

\begin{lem}\label{DT for A to R}
	There is an isomorphism	$Rc_*\bbQ_{A} \simeq \scrI_R^+ \oplus \scrI_R^-$.	Moreover, the intersection complexes $\scrIC_R(\Lambda^1_{R^0} \otimes \scrL_{R^0})$ and $\scrIC_R(\Lambda^3_{R^0} \otimes \scrL_{R^0})$ are sheaves.
	\begin{proof}
		We first discuss some local systems that will appear in the proof. Let $U \subset \Sigma$ be the smooth part of $\Sigma$ and let $\scrL_U$ be the rank one local system on $U$ corresponding to the representation of $\pi_1(U)$ associated to the \'etale double cover $A_U \to U$. Recall that $R^0 = R \setminus NR$ and note that $\Sigma_{R^0} \subset U$. There is an isomorphism of fundamental groups $\pi_1(\Sigma_{R^0}) \simeq \pi_1(U)$ since the complement of $\Sigma_{R^0}$ in $U$ has complex codimension $2$. Due to this isomorphism, we will not distinguish $\scrL_U$ from its restriction to $\Sigma_{R^0}$. There is a surjection of fundamental groups $\pi_1(\Sigma_{R^0}) \to \pi_1(R_0)$  which implies that there is an isomorphism of local systems $\scrL_{U} \simeq r^*\scrL_{R^0}$. 
		
		We now study the Decomposition Theorem for $c \colon A \to R$. Proper base change and the fact that $q$ is a finite map imply that $Rc_* \bbQ_A = Rr_*\bbQ_{\Sigma} \oplus Rr_*(j_{U*}\scrL_U)$. Lemma \ref{DT for A/+-1 to R} implies that there is an isomorphism $Rr_*\bbQ_{\Sigma} \simeq \scrI_{R}^+$. To finish the proof, we must show that there is an isomorphism $Rr_*(j_{U*}\scrL_U) \simeq \scrI_{R}^-$. To see this, notice that by proper base change for $r$ combined with the fact that $j_U$ is an open embedding implies that there is an isomorphism $Rr_*(j_{U*}\scrL_U)|_{R^0} = Rr_*\scrL_{U}$. The isomorphism $\scrL_{U} \simeq r^*\scrL_{R^0}$, coupled with the projection formula and Deligne's Theorem (see Theorem 1.5.3 in the survey \cite{dC-PCMI}), i.e. the Decomposition Theorem for smooth proper morphisms, gives the desired isomorphism $Rr_*(j_{U*}\scrL_U) \simeq \scrI_{R}^-$. The description of the fibers of $c$ imply that there are no additional summands appearing in the Decomposition Theorem for $c$. This description, combined with Relative Hard Lefschetz, also shows that $\scrIC_R(\Lambda^1_{R^0} \otimes \scrL_{R^0})$ and $\scrIC_R(\Lambda^3_{R^0} \otimes \scrL_{R^0})$ are sheaves.
	\end{proof}	
\end{lem}

\subsection{The Relevant Strings over \texorpdfstring{$\Npq$}{Npq}}\label{strings over Npq section}
Recall the Kummer $K3$ diagram from Equation \ref{kummer k3 diagram}. Throughout this section, we work with a fixed $p,q \in J[2]$ with $p \neq q$. We abuse notation and use the same symbols to denote the images of these points in $K$. Let $H$ be a hyperplane section of $K$ and let $E_p$, $E_q$ be the exceptional divisors in the Kummer $K3$ surface $\tildeK$ over nodes $p,q \in K$. Consider the linear system $B' = |\tau^*H - E_p - E_q|$ and note that $B'$ is a base point free pencil. The linear system $B'$ induces an elliptic fibration $e \colon \tildeK \to B'$ with six $I_2$ fibers and two $I_0^*$ fibers (in Kodaira's notation of singular fibers for elliptic surfaces) and admits a section (cf. \cite[\S 4]{Kumar}). We note that the base of this elliptic fibration can be identified with the locus $\Npq$ parameterizing curves in the linear system $B=|2\theta|$ with nodes at $p,q \in J[2]$ via the map $\Npq \to B'$ sending $C$ to $\tau^*f(C) - E_p - E_q$. Under the identification $\Npq \simeq B'$, the six $N(3)$ points in $\Npq$ are identified with the six points in $B'$ having $I_2$ fibers and the two $NR$ points in $\Npq$ are identified with the two points in $B'$ having $I_0^*$ fibers. 

We are now ready to introduce the relevant strings over the loci $\Npq$. Let $\tildecalD_{B'} \subset \tildeK \times B'$ be the universal family of curves associated to the linear system $B'$. The existence of a section induces an identification between the elliptic fibration $e \colon \tildeK_{\Vpq} \to \Vpq$ and the relative Picard scheme $\Pic^0_{\tildecalD_{\Vpq}/\Vpq}$ over the open subset $\Vpq \subset \Npq$. In view of Proposition \ref{irred components over integral curves prop}, this identification implies that there is a surjective morphism of smooth commutative group schemes $G^0|_{\Vpq} \twoheadrightarrow  \tildeK_{\Vpq}$, where $G^0|_{\Vpq}$ is the identity component of $G|_{\Vpq}$, which fiber-by-fiber realizes the Chevalley devissage. We introduce the following complexes, viewable in $D^bMHM_{alg}(\Npq)$ or $D^b(\Npq,\bbQ)$, which we will again call strings. Let $\Lambda^1_{\Vpq} \coloneqq R^1e_*\bbQ_{\tildeK_{\Vpq}}$,
\begin{equation}
\scrI^+_{\Npq} \coloneqq \bbQ_{\Npq} \oplus j_{pq*}(\Lambda_{\Vpq}^1)[-1] \oplus \bbQ_{\Npq}[-2](-1),
\end{equation}
\begin{equation}\label{scrLpq equation}
\scrI^-_{\Npq} \coloneqq j_{pq*} \scrL_{\Vpq} \oplus j_{pq*}(\Lambda_{\Vpq}^1 \otimes \scrL_{\Vpq})[-1] \oplus j_{pq*} \scrL_{\Vpq}[-2](-1),
\end{equation}
where  $\scrL_{\Vpq}$ is a rank one local system on $\Vpq$ with $-1$ monodromy around $NR$ points and trivial monodromy around $N(3)$ points and $j_{pq} \colon \Vpq \to \Npq$ is the inclusion.

\begin{lem}\label{DT for elliptic fibration lem}
	There is an isomorphism
	\begin{equation*}
	Re_*\bbQ_{\tildeK} \simeq \scrI_{\Npq}^+ \oplus (\bbQ_{\Npq \cap N(3)} \oplus \bbQ_{\Npq \cap NR}^{\oplus 4})[-2](-1). 
	\end{equation*}
	\begin{proof}
		This follows immediately from \cite[Example 1.8.4]{dCM1}.
	\end{proof}
\end{lem}

\begin{cor}\label{hodge structure of I_Npq+ corollary}
	There is an isomorphism 
	\begin{equation*}
	H^*(\Npq, \scrI_{\Npq}^+) \simeq H^{ev}(J) \oplus \bbQ^{\oplus 2}[-2](-1)
	\end{equation*} 
	of rational polarizable graded pure Hodge structures.
	\begin{proof}
		Lemma \ref{DT for elliptic fibration lem} gives $H^*(\tildeK) = H^*(\Npq, \scrI_{\Npq}^+) \oplus \bbQ^{\oplus 14}[-2](-1)$. Since $\tildeK$ is the Kummer $K3$ surface associated to $J$, the Hodge structure of $\tildeK$ is given by $H^*(\tildeK) = H^{ev}(J) \oplus \bbQ^{\oplus 16}[-2](-1)$. The result then follows by combining these two equations. 
	\end{proof}
\end{cor}

To understand the string $\scrI_{\Npq}^-$,  we first need to discuss the existence of a special double cover of the Kummer $K3$ surface $\tildeK$, described by Mehran in \cite{Meh}. In each singular fiber of type $I_0^*$, there are four exceptional curves which appear with multiplicity one. Denote these by $E_1, \cdots, E_8$. Mehran shows that $\Delta \coloneqq E_1 + \cdots + E_8 \in 2NS(\tildeK)$. In particular, this means that we can find a double cover $q \colon Z \to \tildeK$ branched along $\Delta$. The preimages $F_i \coloneqq q\inv(E_i)$  are exceptional curves and we can blow them down to obtain a surface $\tau \colon Z \to X$. By \cite[Proposition 2.3]{Meh}. the surface $X$ is the Kummer $K3$ surface associated to some Abelian surface $J'$ which admits an isogeny $J' \to J$ of degree $2$. 

\begin{lem}\label{hodge structure of Z lem}
	There is an isomorphism $H^*(Z) \simeq H^{ev}(J) \oplus \bbQ[-2](-1)^{\oplus 24}$ of rational polarizable graded pure Hodge structures.
	\begin{proof}
		Since $Z$ is the blow up of $X$ at eight points, there is an isomorphism of rational Hodge structures $H^*(Z) \simeq H^*(X) \oplus \bbQ[-2](-1)^{\oplus 8}$. Since $X$ is the Kummer $K3$ associated to the Abelian surface $J'$, there is an isomorphism of rational Hodge structures $H^*(X) \simeq H^{ev}(J')\oplus \bbQ[-2](-1)^{\oplus 16}$. The claim then follows by noticing that the rational Hodge structures of $J'$ and $J$ are isomorphic since $J'$ is isogenous to $J$.
	\end{proof}
\end{lem}

We are now ready to study the string $\scrI_{\Npq}^-$. Let $h \coloneqq e \circ q \colon Z \to \Npq$ denote the composition of the double cover with the elliptic fibration. The number of irreducible components of the fibers of $h$ over $N(3)$ and $NR$ points are straightforward to describe. The fiber of $e$ over an $N(3)$ point is disjoint from the branch locus $\delta$ and is a curve of type $I_2$, i.e. consists of two rational curves which meet transversely at two distinct points. Since $q$ is \'etale over this fiber, each rational component must split into two disjoint rational components. It follows that a fiber of $h$ over an $N(3)$ point has four irreducible components. The fiber of $e$ over a $NR$ is a curve of type $I_0^*$ with five irreducible components. There is a non-reduced component appearing with multiplicity two and four reduced components which lie in the branch locus $\delta$. The non-reduced component, viewed with its reduced induced structure, is a rational curve which intersects the branch locus in four points. It follows that preimage of this rational curve under $q$ is an elliptic curve. In particular, it is connected and it follows that the fiber of $h$ over an $NR$ point consists of five irreducible components. 

\begin{lem}\label{DT for Z}
	There is an isomorphism
	\begin{equation*}
	Rh_*\bbQ_Z \simeq \scrI_{\Npq}^+ \oplus \scrI_{\Npq}^- \oplus \big(\bbQ_{NR}^{\oplus 4} \oplus \bbQ_{N(3)}^{\oplus 2}\big)[-2](-1).
	\end{equation*}
	
	\begin{proof}
		We first discuss some local systems that will appear in the proof. Let $W = \tildeK \setminus \Delta$ and let $\scrL_W$ be the rank one local system on $W$ associated to the \'etale double cover $Z_W \to W$. Recall that $\Upq \coloneqq \Npq \setminus NR \simeq \bbC^*$ and note that the inclusion $\tildeK_{\Upq} \subset W$ holds since $\Delta$ is contained within the two $I_0^*$ fibers of $e$. Let $\scrL_{\tildeK_{\Upq}}$ be the restriction of $\scrL_W$ to $\tildeK_{\Upq}$ and let $\scrL_{\Upq}$ be the rank one local system corresponding to the representation $\pi_1(\Upq) \to \Aut(\bbQ)$ where the generator acts by $-1$. One can show that $\scrL_{\tildeK_{\Upq}} \simeq e^*\scrL_{\Upq}$ using the surjection $\pi_1(\tildeK_{\Upq}) \twoheadrightarrow \pi_1(\Upq)$ between fundamental groups. 
		
		We now study the Decomposition Theorem for $h \colon Z \to \Npq$. Proper base change and the fact that $q$ is a finite map imply that $Rh_*\bbQ_Z = Re_*\bbQ_{\tildeK} \oplus Re_*(j_{W*} \scrL_W)$. Lemma \ref{DT for elliptic fibration lem} implies that there is an isomorphism $Re_*\bbQ_{\tildeK} \simeq \scrI_{\Npq}^+$. To finish the proof, note that proper base change for $e$ implies there is an isomorphism $j_{pq}^*Re_*(j_{W*} \scrL_W) \simeq Re_*\scrL_{\tildeK_{\Upq}}$. This isomorphism, coupled with the projection formula and Deligne's Theorem (see Theorem 1.5.3 in the survey \cite{dC-PCMI}), imply that the strings $\scrI_{\Npq}^+$ and $\scrI_{\Npq}^-$ are direct summands appearing in the Decomposition Theorem for $Rh_*\bbQ$. Considering the description of fibers of $h$ over $N(3)$ and $NR$ points above gives the skyscraper contributions.
	\end{proof}
\end{lem}

\begin{cor}\label{hodge structure of I_Npq- corollary}
	There is an isomorphism $H^*(\Npq, \scrI_{\Npq}^-) \simeq \bbQ^{\oplus 2}[-2](-1)$ of rational polarizable graded pure Hodge structures.
	\begin{proof}
	 We have $H^*(Z) \simeq H^{ev}(J) \oplus H^*(\Npq, \scrI_{\Npq}^-) \oplus \bbQ[-2](-1)^{\oplus 22}$ by Lemma \ref{DT for Z} and Corollary \ref{hodge structure of I_Npq+ corollary}. The result then follows from Lemma \ref{hodge structure of Z lem},
	\end{proof}
\end{cor}

\subsection{Ng\^o Strings for the Lagrangian Fibrations \texorpdfstring{$\tildem$}{tildem}, \texorpdfstring{$n$}{n}}\label{DT for tildeM, N section}

\begin{prop}\label{DT for tildeM and N}
	Let $\langle \bullet \rangle \coloneqq [-2\bullet](-\bullet)$. The Decomposition Theorems for $R\widetilde{m}_*\bbQ_{\tildeM}$, $Rn_*\bbQ_N$ in $D^bMHM_{alg}(B)$ take the following form: 
	\begin{equation}\label{DT for tildeM equation}
		R\widetilde{m}_*\bbQ_{\tildeM} \simeq \scrI_{B} \oplus \scrI_{R}^+\langle 1 \rangle \oplus \bigoplus_{p,q \in J[2], p \neq q} \scrI_{\Npq}^+\langle 2 \rangle \oplus \bbQ_{NR}^{\oplus r + 16}\langle 3 \rangle
	\end{equation}
	
	\begin{equation}\label{DT for N equation}
		Rn_*\bbQ_{N} \simeq \scrI_{B} \oplus \scrI_{R}^-\langle 1 \rangle \oplus \bigoplus_{p,q \in J[2], p \neq q} \scrI_{\Npq}^-\langle 2 \rangle \oplus \bbQ_{NR}^{\oplus r}\langle 3 \rangle.
	\end{equation}	
	where $r = 0$ or $1$.
	\begin{proof}
		Since the triples $(\tildeM, B, G)$ and $(N, B, G)$ are $\delta$-regular weak Abelian fibrations by Proposition \ref{(tildeM,B,G), (N,B,G) are drwaf}, we can apply the Ng\^o Support Theorem.
		
		We first deal with $R\widetilde{m*}\bbQ_{\tildeM}$. Proposition \ref{DT for symplectic resolution}, combined with Lemmas \ref{DT for A/+-1 to R} and \ref{type I sheaf}, implies that there is an isomorphism $R\tildem_*\bbQ_{\tildeM} \simeq Rm_*\scrIC_{M} \oplus \scrI_R^+\langle 1 \rangle \oplus \bbQ_{NR}^{\oplus 16} \langle 3 \rangle$. In particular, $B$, $R$, and $NR$ must be supports. Moreover, the Ng\^o Support Theorem implies that the local systems appearing in the string $\scrI_R^+$ can be identified with the cohomology of the Abelian part of the identity component of the group scheme $G$ over the Zariski open subset $R^0 \subset R$. The direct sum of the Ng\^o strings associated with these three supports is:
		\begin{equation}\label{direct sum of Ngo strings for tildeM associated to supports S,R,NR}
			\scrI_{B}^{\oplus r_{\tildeM,B}} \oplus {\scrI_{R}^+}^{ \oplus r_{\tildeM,R}}\langle 1 \rangle \oplus \bbQ_{NR}^{r_{\tildeM,NR}}\langle 3 \rangle.
		\end{equation}
		for some strictly positive integers $r_{\tildeM, B}, r_{\tildeM, R}$, and $r_{\tildeM,NR}$. Using Lemmas \ref{cohomology sheaf of I_B} and \ref{DT for A to R} to compute the combined contribution of these summands to the highest direct image $R^6_{\tildeM}$, Proposition \ref{direct image sheaves in top degree prop} implies that $r_{\tildeM, |2\theta|} = r_{\tildeM, R} = 1$ and $16 \leq r_{\tildeM, NR} \leq 34$.

		For any $p,q \in J[2]$ with $p \neq q$, the contribution of the summands in Equation \ref{direct sum of Ngo strings for tildeM associated to supports S,R,NR} to the highest direct image $R^6_{\tildeM}$ restricted to the locus $\Vpq$ is simply $\bbQ_{\Vpq}$. Since  $R^6_{\tildeM|_{\Vpq}} \simeq \bbQ_{\Vpq}^{\oplus 2}$ by Proposition \ref{direct image sheaves in top degree prop}, $\Npq$ must be a support and the associated Ng\^o string is $\scrI_{\Npq}^+\langle 2 \rangle$. Thus, we have shown that the direct sum of the Ng\^o strings appearing in the Decomposition Theorem for $\tildeM$ associated with the supports $B$, $R$, $\Npq$, and $NR$ is
		\begin{equation}\label{direct sum of Ngo strings for tildeM associated to supports S,R,Npq,NR}
			\scrI_{B} \oplus \scrI_{R}^+\langle 1 \rangle \oplus \bigoplus_{p,q \in J[2], p \neq q} \scrI_{\Npq}^+\langle 2 \rangle \oplus \bbQ_{NR}^{\oplus r_{\tildeM,NR}}\langle 3 \rangle
		\end{equation}
		
		According to our description of $R^6_{\tildeM}$, the only other possible support is $N(3)$. If $C \in N(3)$ has nodes at distinct points $p,q,r \in J[2]$, then $C =\Npq \cap N_{pr} \cap N_{qr}$. It follows that the restriction of Equation \ref{direct sum of Ngo strings for tildeM associated to supports S,R,Npq,NR} to $C$ is $\bbQ^{\oplus 4}$. Since $R^6_{\tildeM}|_C \simeq \bbQ^{\oplus 4}$ by Proposition \ref{direct image sheaves in top degree prop}, $C$ is not a support and Equation \ref{DT for tildeM equation} holds.
		
		We note that the skyscraper sheaf $\calH^2(\scrIC_B(\Lambda_B^4))$ can thus only be supported on $NR$. The stalks of this sheaf over different points of $NR$ are all isomorphic, so there exists an integer $r_{2,4} \geq 0$ such that $\calH^2(\scrIC_B(\Lambda_B^4)) = \bbQ_{NR}^{\oplus r_{2,4}}$. Although the presence of this skyscraper sheaf prohibits us from determining the exact shape of the Decomposition Theorem, we can say the following. Since any $C \in NR$ passes through exactly six 2-torsion points, $C$ lies in $15 = \binom{6}{2}$ lines of the form $\Npq$. Restricting Equation \ref{direct sum of Ngo strings for tildeM associated to supports S,R,Npq,NR} to $C$ then gives $r_{\tildeM, NR} = 17 - r_{2,4}$. 
		
		We now deal with $Rn_*\bbQ_{N}$. Proposition \ref{direct image sheaves in top degree prop} implies that $B$ must be a support and the associated Ng\^o string is $\scrI_B$. Lemma \ref{cohomology sheaf of I_B} and the fact that the skyscraper sheaf $\calH^2(\scrIC_B(\Lambda_B^4))$ is supported on $NR$ imply that this contribution restricted to $R^0$ is $\bbQ_{R^0}$. Since ${R^6_N}|_{R^0} \simeq \bbQ_{R^0} \oplus \scrL_{R^0}$ by Proposition \ref{direct image sheaves in top degree prop}, we see that $R$ is must be a support and the associated Ng\^o string is $\scrI_{R}^-$.
		
		By Lemmas \ref{cohomology sheaf of I_B} and \ref{DT for A/+-1 to R}, for any $p,q \in J[2]$ with $p \neq q$, the contribution of the strings $\scrI_B$ and $\scrI_{R}^-$ restricted to $\Vpq$ is $\bbQ_{\Vpq}$. Since ${R^6_N}|_{\Vpq} \simeq \bbQ_{\Vpq} \oplus \scrL_{\Vpq}$ by Proposition \ref{direct image sheaves in top degree prop}, we see that $\Npq$ must be a support and the associated Ng\^o string is $\scrI_{\Npq}^-$. We have thus shown that the direct sum of the Ng\^o strings appearing in the Decomposition Theorem for $N$ associated with the supports $B$, $R$, and $\Npq$ is
		\begin{equation}\label{direct sum of Ngo strings for N associated to supports S,R,Npq}
			\scrI_{B} \oplus \scrI_{R}^-\langle 1 \rangle \oplus \bigoplus_{p,q \in J[2], p \neq q} \scrI_{\Npq}^-\langle 2 \rangle.
		\end{equation}
		
		According to our description of $R^6_{N}$, the only other possible supports are $N(3)$ and $NR$. Using the same argument as for $\tildeM$, one can see that $N(3)$ is not a support. If $C \in NR$, then the description of the local systems $\scrL_{R^0}$ and $\scrL_{\Vpq}$ imply that restriction of Equation \ref{direct sum of Ngo strings for N associated to supports S,R,Npq} to $C$ is $\bbQ$. Since $R^6_{N}|_C \simeq \bbQ^{\oplus 2}$ by Proposition \ref{direct image sheaves in top degree prop}, $C$ is potential support. Noting that $\calH^2(\scrIC_B(\Lambda_B^4)) = \bbQ_{NR}^{\oplus r_{2,4}}$, we conclude that $r_{N,NR} = 1-r_{2,4}$. Letting $r = r_{N,NR} = 1-r_{2,4}$, we see that $r_{\tildeM,NR} = r+16$ as desired. We note that $r=0$ or $1$ depending on the contribution $r_{2,4}$ from $\calH^2(\scrIC_B(\Lambda_B^4))$. 
	\end{proof}
\end{prop}

\subsection{The Cohomology of \texorpdfstring{$\tildeM$}{tildeM}}\label{proof of numerical main thm}
\begin{rmk}\label{cat of HS rmk}
	Let $\calA$ be a semisimple Abelian category where every object has finite length and the isomorphism classes of simple objects form a set $\frakS$. Every object $a \in \calA$ is isomorphic to a unique finite direct sum $a \simeq \bigoplus_{\fraks \in \frakS} \fraks^{\oplus n_{\fraks}(a)}$ of simple objects with multiplicities $n_{\fraks}(a)$.	If we have an identity $[a] = [b]-[c]$ in the Grothendieck group $K(\calA)$ with $a,b,c \in \calA$, then $n_{\fraks}(a) = n_{\fraks}(b) - n_{\fraks}(c)$.	Let $\phi \colon \Obj (\calA) \to \frakM$ be an assignment into a commutative group which is additive in exact sequences. If $[a] = [b]-[c]$ as above, then $\phi(a) = \phi(b) - \phi(c)$.
\end{rmk}

\begin{prop}\label{coh of tildeM}
	Recall that $A$ is the Abelian four-fold $A \coloneqq J^\vee \times J$. Let $\langle \bullet \rangle \coloneqq [-2\bullet](-\bullet)$. Then we have
	\begin{equation*}
	H^*(\tildeM) = H^*(N) + \big(H^{ev}(A)^{\oplus 2} - H^*(A) \big)\langle 1 \rangle + H^{ev}(J)^{\oplus 120}\langle 2 \rangle + H^*(NR)^{\oplus 16}\langle 3 \rangle
	\end{equation*}
	in the Grothendieck group of rational graded polarizable pure Hodge structures.
	\begin{proof}
		Noting that $H^{ev}(A) = H^*(A/\pm1)$ and $H^{ev}(J) = H^*(J/\pm1)$, this follows from Proposition \ref{DT for tildeM and N}, Lemmas \ref{DT for A/+-1 to R} and \ref{DT for A to R}, and Corollaries \ref{hodge structure of I_Npq+ corollary} and \ref{hodge structure of I_Npq- corollary}. 
	\end{proof}
\end{prop}

\begin{cor}
	We recover the known Betti and Hodge numbers for manifolds of $OG6$ type, which were first computed in \cite{MRS1}.
	\begin{proof}
		Noting that the Betti and Hodge numbers of $N$ can be computed from the G\"ottsche-Soergel formula in \cite{GS}, this follows from Proposition \ref{coh of tildeM} and Remark \ref{cat of HS rmk}. 
	\end{proof}
\end{cor}

\begin{thm}\label{main thm hodge structure}
	Let $U = H^{ev}(J,\bbQ)$ and $W = H^{odd}(J, \bbQ)$. Then
	\begin{equation*}
	H^*(\tildeM) = \Sym^3U \oplus \big( (U^{\otimes 2})^{\oplus 2} \oplus W^{\otimes 2} \big) \langle 1 \rangle  \oplus U^{\oplus 137}\langle 2 \rangle \oplus \bbQ^{\oplus 512}\langle 3 \rangle.
	\end{equation*}
	
	\begin{proof}
		The description of the rational Hodge structure of $N$ in terms of the rational Hodge structure of $J$ can be read off from the $LLV$ decomposition of the cohomology of $N$ given in Corollary 3.6 of \cite{Green-Kim-Laza-Robles}. Using the branching rules described in \cite[\S 25.3]{FH}, one computes that the rational Hodge structure of $N$ is given by
		\begin{equation*}
		H^*(N) = \Sym^3U \oplus U^{\otimes 2}\langle 1 \rangle \oplus U^{\oplus 16}\langle 2 \rangle \oplus \bbQ^{\oplus 256}\langle 3 \rangle \oplus (U\otimes W)^{\oplus 2}\langle 1 \rangle. 
		\end{equation*}
		
		Using the identification of the rational Hodge structure $H^*(J^\vee)$ with $H^*(J)$ via the principal polarization $\theta$, the rational Hodge structure of $A = J^\vee \times J$ is $H^*(A) = U^{\otimes 2} \oplus \big(U \otimes W)^{\oplus 2} \oplus W^{\otimes 2}$. The rational Hodge structure of $A / \pm 1$ is the $\bbZ/2\bbZ$-invariant part of $H^*(A)$, which is given by $H^*(A/ \pm 1) = U^{\otimes 2} \oplus W^{\otimes 2}$. Proposition \ref{coh of tildeM} then gives the result.
	\end{proof}
\end{thm}

\begin{rmk}\label{cohomology of tildeM in terms of even cohomology of J}
	Using Schur functors (see \cite[Ch. 6]{FH}), one can express the Hodge structure of $\tildeM$ purely in terms of $U$. In particular, one can show that 
	\begin{equation*}
		H^*(\tildeM) = \Sym^3U \oplus \Lambda^3U \oplus (U^{\otimes 2})^{\oplus 2}\langle 1 \rangle \oplus U^{\oplus 138}\langle 2 \rangle \oplus \bbQ^{\oplus 512}\langle 3 \rangle,
	\end{equation*}
	which agrees with the description given in \cite{Green-Kim-Laza-Robles}. 
\end{rmk}

\bibliographystyle{halpha} 
\small{\bibliography{OG6bib}}

\end{document}